\documentclass[11pt]{article}
\usepackage{amsfonts, amsmath, amssymb, graphics, epsfig}
\usepackage{xcolor}
\usepackage{enumerate, amstext, amsthm,bm}

\theoremstyle{plain}
\newtheorem{theorem}{Theorem}[section]
\newtheorem{lemma}[theorem]{Lemma}
\newtheorem{corollary}[theorem]{Corollary}
\newtheorem{proposition}[theorem]{Proposition}
\newtheorem{remark}[theorem]{Remark}

\theoremstyle{definition}

\newtheorem{example}[theorem]{Example}
\theoremstyle{remark}

\def\IC{{\mathbb C}}

\def\IN{{\mathbb N}}
\def\IR{{\mathbb R}}
\def\bA{{\mathbf A}}
\def\bB{{\mathbf B}}
\def\bC{{\mathbf C}}
\def\bD{{\mathbf D}}

\def\bQ{{\mathbf Q}}
\def\bK{{\mathbf K}}
\def\bM{{\mathbf M}}
\def\bS{{\mathbf S}}
\def\bT{{\mathbf T}}

\def\bV{{\mathbf V}}

\def\bu{{\mathbf u}}
\def\bv{{\mathbf v}}

\def\bd{{\mathbf d}}
\def\b0{{\mathbf 0}}
\def\bw{{\mathbf w}}
\def\bx{{\mathbf x}}

\def\bmu{{\bm{\mu}}}
\def\bxi{{\bm{\xi}}}

\def\cB{{\mathcal B}}

\def\cQ{{\mathcal Q}}
\def\cD{{\mathcal D}}

\def\cW{{\mathcal W}}

\def\cK{{\mathcal K}}

\def\cP{{\mathcal P}}
\def\cS{{\mathcal S}}

\def\cH{{\mathcal H}}
\def\SH{{\cS(\cH)}}

\def\KH{{\cK(\cH)}}
\def\BH{{{\mathcal B}(\cH)}}
\def\re{{\rm Re}\,}

\def\tr{{\rm tr}\,}

\def\span{{\rm span}\,}

\def\diag{{\rm diag}\,}
\def\conv{{\bf conv}\,}
\def\cl{{\bf cl}\,}

\def\[{\left [}
\def\]{\right ]}
\def\({\left (}
\def\){\right )}
\def\ess{{\rm ess}}
\def\la{{\langle}}
\def\ra{{\rangle}}
\def\bmu{{\mathbf \mu}}
\def\bR{{\mathbf R}}

\def\<{{\langle}}
\def\>{{\rangle}}
\def\1{{\mathbf 1}}

\def\bnu{{\bm \nu}}
\def\bmu{{\bm \mu}}
\def\bxi{{\bm \xi}}
\def\b0{{\bf 0}}

\def\ess{{\rm ess}}

\begin{document}
\openup .7\jot

\title{Joint $k$-numerical ranges of operators}

\author{Jor-Ting Chan, Chi-Kwong Li, Yiu-Tung Poon}

\date{}
\maketitle
\begin{abstract}
Let ${\mathcal H}$ be a complex Hilbert space
and   ${\mathcal B}({\mathcal H})$  the algebra of
all bounded linear operators on ${\mathcal H}$.
For a positive integer $k$ less than the dimension of
${\mathcal H}$ and ${\mathbf A} = (A_1, \dots, A_m)\in {\mathcal B}({\mathcal H})^m$,
the joint $k$-numerical range $W_k({\mathbf A})$ is the
set of $(\alpha_1, \dots, \alpha_m) \in{\mathbb C}^m$ such that
$\alpha_i = \sum_{j = 1}^k \langle A_ix_j, x_j\rangle$
for an orthonormal set $\{x_1, \ldots,  x_k\}$ in ${\mathcal H}$.
The geometrical properties of $W_k(\bA)$ and their
relations with the algebraic properties of $\{A_1, \dots, A_m\}$
are investigated in this paper.
First, conditions for $W_k({\mathbf A})$ to be convex are studied, and some results
for finite dimensional operators are extended to the infinite dimensional
case. An example of $\bA$ is constructed such that $W_k({\mathbf A})$
is not convex, but $W_r({\mathbf A})$ is convex for all  positive integer $r$
not equal to $k$.
Second, descriptions are given for the closure of $W_k({\mathbf A})$ and
the closure of $\conv W_k({\mathbf A})$ in terms of the joint essential
numerical range of ${\mathbf A}$ for infinite dimensional operators $A_1, \dots, A_m$.
These lead to characterizations of $W_k({\mathbf A})$
or $\conv W_k({\mathbf A})$ to be closed. Moreover, it is shown that
$\conv W_k({\mathbf A})$ is closed whenever $W_{k+1}(\bA)$ or
${\rm conv} W_{k+1}({\mathbf A})$ is.
These results are used to study the connection between the geometric
properties of $W_k({\mathbf A})$ and algebraic  properties of $A_1, \dots, A_m$.
For instance,   $W_k({\mathbf A})$ is a polyhedral set,
i.e., the convex hull of a finite set, if and only if
$A_1, \dots, A_k$ have a common reducing subspace ${\mathbf V}$
of finite dimension such that the compression of $A_1, \dots, A_m$ on the subspace
${\mathbf V}$ are diagonal operators $D_1, \dots, D_m$ and
$W_k({\mathbf A}) = W_k(D_1, \dots, D_m)$. Characterization is also
given to ${\bf A}$ such that the closure of $W_k({\mathbf A})$ is polyhedral.
For finite rank operators the following two condition are
equivalent:
(1) $\{A_1, \dots, A_m\}$ is a commuting family of normal operators.
(2) $W_k(A_1, \dots, A_m)$ is polyhedral for every positive integer $k$ less
than $\dim \cH$.
However, the two conditions are not equivalent for compact operators.
Characterizations are given
for compact operators $A_1, \dots, A_m$ satisfying (1) and (2), respectively.
Results are also obtained for general non-compact operators.
Open problems and future research topics are presented.
\end{abstract}



\section{Introduction}\label{S:1}

Let $\mathcal H$ be a complex Hilbert space
equipped with the inner product $\la x, y\ra$,  and let $\BH$
be the  algebra of all bounded linear operators
on $\mathcal H$.   If $\dim \cH = n$ is finite, we identify
$\cH$ with $\IC^n$, $\BH$ with $M_n$, and $\la x, y\ra = y^*x$.
The numerical range of $A \in \BH$ is defined as
$$W(A) = \{\langle A x, x\rangle: \mbox{$x$ is a unit vector in
$\cH$}\},$$
which is a useful tool for studying operators
and matrices; see \cite{H,HJ} and its references.
For instance, for $A \in \BH$, $W(A) = \{\mu\}$ if and only if $A = \mu I$;
$A \subseteq \IR$ if and only if $A = A^*$; if $A$ is normal then
$ \cl(W(A)) = \conv\sigma(A)$, the convex hull of the spectrum $\sigma(A)$ of $A$.
Here and in the sequel,
we shall write $\cl X$ and $\conv X$
to denote the closure and convex hull of $X$ respectively, for any subset $X$ of a normed space.

Motivated by theory and applications, researchers consider different
generalizations of the classical numerical range.
For instance, for a positive integer $k \le \dim \cH$, Halmos \cite{H}
introduced the $k$-numerical range of  $A \in \BH$ defined by
$$W_k(A) = \left\{ \sum_{j=1}^k \la A x_j, x_j\ra:
\{x_1,\dots, x_k\}\subseteq \cH \hbox{ is an orthonormal set}\right\}.$$
Evidently,  $W_1(A)$ reduces to the classical numerical range $W(A)$.
Researchers have extended the results on $W(A)$ to $W_k(A)$, and showed
that $W_k(A)$ can provide additional information about the operator;
see later sections, \cite{AT,BFL,B,BW,H,HJ,LP1,LST,MF,MMF} and their references.
For instance, it is known that if $A \in M_n$ is normal, then $W(A)$ is a
polygon equal to $\conv\sigma(A)$. But the converse is not true.
Nevertheless, it was shown in \cite{LST} that $A\in M_n$ is normal
if and only if $W_k(A)$ is a polygon for some $k$ with $|n - k/2| \le 1$.

In the study of multiple operators
$A_1, \dots, A_m$, researchers consider the joint $k$-numerical range of
$\bA = (A_1, \dots, A_m)\in \BH^m$ defined and denoted by
$$W_k(\bA)=\left\{
\left( \sum_{\ell = 1}^k \langle A_1x_\ell, x_\ell\rangle,
\dots, \sum_{\ell = 1}^k \langle A_mx_\ell, x_\ell\rangle\right):
\{x_1, \dots, x_k\} \subseteq \cH \hbox{ is an orthonormal set} \right\}$$
for $k \le \dim \cH$. When $k = 1$, we use the notation
$W(\bA)$. For simplicity, we often write $\la\bA x, x\ra = (\la A_1 x, x\ra,
\dots, \la A_m x, x\ra)$ for $x\in \cH$. Then elements of
$W_k(\bA)$ are of the form $\sum_{\ell = 1}^k \la \bA x_\ell, x_\ell\ra$ for an
orthonormal set $\{x_1, \dots, x_k\}$ in $\cH$.
One may see some background for $W(\bA)$ in \cite{AT,BFL,BW,LP1,LPW,Juneja,Muller}
and their references. Useful
information about $\{A_1,\dots, A_m\}$ can be obtained from  $W_k(\bA)$.
For instance, it was shown in \cite{LPW} that
$\{A_1, \dots, A_m\}\subseteq M_n$ is a commuting family
of normal matrices if and only if there is $k$ with $|n/2-k|\le 1$
such that $W_k(A_1, \dots, A_m))$ is a polyhedral set, i.e., the
convex hull of a finite set in $\IC^m$.

In this paper, we obtain additional results
connecting the geometrical properties of $W_k(\bA)$ and the
algebraic properties of $\{A_1, \dots, A_m\}$ as well as showing that
direct extensions of some
results on $W_k(A)$ for a single operator $A$ to $W_k(\bA)$ is impossible.

We begin by collecting some basic results in Section 2.
In Section 3, we consider the convexity of $W_k(\bA)$.
The celebrated Toeplitz-Hausdorff theorem asserts that
$W(A)$ is always convex.  Berger extended the result to
$W_k(A)$ for any $k \le \dim \cH$; for example, see \cite{H}.
However, for $\bA\in \BH^m$, $W_k(\bA)$ may not be convex if $m \ge 2$.
We will extend
some results concerning the convexity of $W_k(\bA)$
for finite dimensional operators to the infinite dimensional case.
An example of $\bA$ is constructed such that
$W_k(\bA)$ is not convex, but $W_r(\bA)$ is convex  for all  positive integer $r$
not equal to $k$. In particular, this shows that there are no implications between
the convexity of  $W_k(\bA)$ and $W_{k+1}(\bA)$.
This is quite different from the other properties considered in the subsequent sections,
namely, special properties of $W_{k+1}(\bA)$ often induces similar properties on $W_k(\bA)$.

For infinite dimensional operator, $W_k(\bA)$ may not be closed.
In Section 4, we determine the closure of $W_k(\bA)$
in terms of the joint essential numerical
range of $\bA$ defined by
\begin{equation}
\label{Wess}
W_{\ess}(\bA)= \cap\{\cl(W_1(\bA+\bK)):\bK\in \cK(\cH)^m\},
\end{equation}
where $\cK(\cH)$ is the set of all compact operators in $\BH$.
These lead to characterizations of $W_k(\bA)$ or
$\conv W_k(\bA)$ to be closed. Moreover, it is shown that
$\conv W_k(\bA)$ is closed whenever $W_{k+1}(\bA)$ or $\conv W_{k+1}(\bA)$
is.
The results extend those in \cite{CLP0,L}, and will be useful in
the subsequent sections, which concern
the interplay between the geometric
properties of $W_k(\bA)$ and
algebraic properties of $A_1, \dots, A_m$.

In section 5, we show that $W_k(\bA)$ is a polyhedral set in $\IC^m$
if and only if
$A_1, \dots, A_k$ have a common reducing subspace ${\mathbf V}$
of finite dimension
such that the compression of $A_1, \dots, A_m$ on the subspace
${\mathbf V}$ are diagonal operators $D_1, \dots, D_m$ and
$W_k({\mathbf A}) = W_k(D_1, \dots, D_m)$.
This extends the result in \cite{LPW}.
Characterization is also
given to ${\bf A}$ such that the closure of $W_k({\mathbf A})$ is polyhedral.

For finite rank operators, one can show that the following two conditions are
equivalent:

\medskip
(a) $\{A_1, \dots, A_m\}$ is a commuting family of normal operators.

\medskip
(b) $W_k(A_1, \dots, A_m)$ is polyhedral for every positive integer $k < \dim \cH$.

\medskip\noindent
This result is no longer true for compact operators.  In Section 6, characterizations are given
for compact operators $A_1, \dots, A_m$ satisfying (a) and (b), respectively.
Additional results concerning conditions (a) and (b) are obtained for general
non-compact operators in Section 7.

In Section 8, we describe some related results, open problems
and additional research topics.

In our discussion, we will
use $\BH,\SH,\KH$ to denote the set of bounded linear operators, the set of
self-adjoint operators, and the set of compact operators acting on $\cH$.
Let $\cH_1 $ and $\cH_2$ be Hilbert spaces. A linear map $X:\cH_1 \to \cH_2$
is an isometry if $X^*X = I_{\cH_1}$.
Denote by $M_{m,n}$ and $M_n$ the set of $m\times n$ matrices and the set of
$n\times n$ matrices.
We will use the convention that $W_0(\bA) = \{(0,\dots, 0)\}.$
Suppose $\dim \cH = n$ is finite. Then
$W_n(\bA) = \{(\tr A_1, \dots, \tr A_m)\}$, and
$W_k(\bA) = (\tr A_1,\dots, \tr A_m) - W_{n-k}(\bA).$
Unless it is specified otherwise,  we always assume that
$k$ is a positive integer such that $1\le k < \dim \cH$
when $W_k(\bA)$ is considered. We note that $(a_1, \dots, a_m) \in W_k(\bA)$
if and only if any of the following holds.
\begin{itemize}
\item There is a rank $k$ orthogonal projection $P$ such that
$(a_1,\dots, a_m) = (\tr(A_1P), \dots, \tr(A_mP))$.
\item There is an isometry $X: \IC^k \rightarrow \cH$ such that
$(a_1, \dots, a_m) = (\tr(X^*A_1X), \dots, \tr(X^*A_mX))$.
\end{itemize}
Clearly, $X: \IC^k\to \cH$ is an isometry such that
$X^*AX = B$ if and only if there is an orthonormal basis
$\{x_1, \dots, x_k\} \cup \cB_2$
of  $\cH$  such that $B = (\la Ax_j, x_i\ra)$.
In such a case, we will say that there is a unitary $U \in \BH$ such that
$U^*AU = \begin{pmatrix} B & \star \cr \star & \star\end{pmatrix}.$
With this notation, we see that
$(a_1, \dots, a_m) \in W_k(A)$ if and only if the following holds.
\begin{itemize}
\item There is a unitary $U \in \BH$ such that
$
U^*A_jU  = \begin{pmatrix} \tilde A_j & \star \cr \star & \star \cr \end{pmatrix}
$
with $\tilde A_j \in M_k$
such that $\tr \tilde A_j = a_j$ for $j = 1, \dots, m$.
\end{itemize}

In our study, it will be shown that if
$(a_1, \dots, a_m) \in W_k(\bA)$ has some special properties, then
there will be a unitary $U \in \BH$ such that $U^*A_jU = \begin{pmatrix}
\tilde A_j & 0 \cr 0 & \star\cr\end{pmatrix}$. Thus,
$\span\{x_1, \dots, x_k\}$ is a reducing subspace for $A_j$ for $j = 1,\dots, m$.
In such a case, we will say that $A_1, \dots, A_m$ have a common reducing space
$\cH_1$ of $\cH$ of dimension $k$ and $X:\cH_1\rightarrow \cH$ is an isometry
such that $(a_1, \dots, a_m) = (\tr(X^*A_1X), \dots, \tr(X^*A_mX))$.

\section{Basic properties}
\setcounter{equation}{0}


\begin{proposition} \label{bounded}
Let $\bA = (A_1, \dots, A_m) \in \BH^m$. Then $W_k(\bA)$ is bounded.
If $\dim \cH$ is finite, then $W_k(\bA)$ is closed.
\end{proposition}

For any positive integer $\ell$, let
$I_\ell$ and $0_\ell$ be the $\ell\times \ell$ identity and zero matrices respectively.

\begin{example} Let $A = I_{k} \oplus 0_{k} \oplus \diag(1,1/2, 1/3, \dots)$.
It is easy to check that $W_k(A) = [0,k]$ is closed, but
$W_{k+1}(A) = (0, k+1]$ is not.
\end{example}

In Section 4, we will consider the condition for $W_k(\bA)$ to be closed for
$\bA \in \BH^m$.

\medskip
Let $A = H+iG \in \BH$, where $H = (A+A^*)/2$ and $G = (A-A^*)/(2i)$ are self-adjoint.
Then $W(A) = W(H, G)$ if we identified $\IC = \IR^2$.  If
$X: \cH_1 \rightarrow \cH$ is an isometry,
then $W(X^*AX) \subseteq W(A)$. If $U$ is unitary, then $W(A) = W(U^*AU)$.
These properties have direct extensions to $W_k(\bA)$.

\begin{proposition} \label{H+iG}
Let $\bA = (A_1, \dots, A_m) \in \BH^m$.
Suppose $A_j = H_j + iG_j$ for $j = 1, \dots, m$, where
$H_j, G_j$ are self-adjoint. If we identify $\IC^m$ with $\IR^{2m}$, then
$W_k(\bA)$ can be identified with
$$W_k(H_1,G_1, \dots, H_m,G_m) \subseteq \IR^{2m}.$$
\end{proposition}

By the above proposition, to understand the geometrical properties
and topological properties of $W_k(\bA)$, one can focus on self-adjoint
operators $A_1, \dots, A_m$.

\begin{proposition} \label{X*AX}
Let $\bA = (A_1, \dots, A_m) \in \BH^m$.
Suppose $X: \cH_1 \rightarrow \cH$ is an isometry,
 where $\dim \cH_1 \ge k$. Then
$$W_k(X^*A_1X, \dots, X^*A_mX) \subseteq W_k(A_1, \dots, A_m).$$
In particular, if $U \in \BH$ is unitary, then
$$W_k(U^*A_1U, \dots, U^*A_mU) = W_k(\bA).$$
\end{proposition}

Suppose $\alpha, \beta \in \IC$ and $A \in \BH$. Then
$W_k(\alpha A) = \alpha W_k(A)$ and $W_k(A + \beta I) = W_k(A) + k\beta$.
We have the following extension.

\begin{proposition} \label{affine}
Let $\bA = (A_1, \dots, A_m) \in \BH^m$.
\begin{itemize}
\item[{\rm (a)}] If $(\mu_1, \dots, \mu_m)\in \IC^m$,
then
$W_k(A_1+\mu_1 I, \dots, A_m + \mu_m I)$ equals
$$W_k(\bA) + k(\mu_1,\dots, \mu_m)
= \{(a_1, \dots, a_m) +
k(\mu_1, \dots, \mu_m): (a_1, \dots, a_m) \in W_k(\bA)\}.$$
\item[{\rm (b)}] If
 $T = (T_{ij})$ be an $m\times p$ matrix, and
$(B_1, \dots, B_p) =  (A_1, \dots, A_m) (T_{ij}I_\cH)$,
then
$$W_k(B_1, \dots, B_p) = \{(a_1, \dots, a_m)T: (a_1, \dots, a_m) \in W_k(\bA)\}.$$
\item[{\rm (c)}] There are $C_1, \dots, C_q \in \BH^q$ such that
$\{C_1, \dots, C_q,I\}$ is linearly independent, $R = (R_{ij})\in M_{q,m}$ and $(s_1,\dots,s_m)\in \IC^m$
such that
$(A_1, \dots, A_m) = (C_1, \dots, C_q)(R_{ij}I_\cH)+(s_1I_\cH ,\dots,s_mI_\cH)$ and
$$W_k(\bA) = \{ (c_1, \dots, c_q)R: (c_1, \dots, c_q) \in W_k(C_1, \dots, C_q)\} +k(s_1,\dots,s_m).$$
\end{itemize}
\end{proposition}

Several remarks are in order.
\begin{itemize}

\item Proposition \ref{affine}
asserts that
if one applies a translation $\bA \mapsto \bB = \bA + (\mu_1 I, \dots, \mu_m I)$,
then $W(\bB)$ is obtained from $W(\bA)$ by
the translation $(a_1, \dots, a_m) \mapsto (a_1, \dots, a_m) + k(\mu_1, \dots, \mu_m)$;
if one applies a  linear transform
$\bA \mapsto \bB = \bA(T_{ij} I_{\cH})$,
then $W(\bB)$ is obtained from $W_k(\bA)$ by the same transformation
$(a_1,\dots, a_m) \mapsto (a_1,\dots, a_m)T$.
Hence, the study of the geometric properties $W_k(\bA)$ can be reduced to the study of
$W_k(C_1, \dots, C_q)$, where $\{C_1, \dots, C_q, I\}$
is a linearly independent set.  These transformations will be used
frequently in our subsequent discussion.
\item
In Proposition \ref{affine}, suppose
$A_1, \dots, A_m$ are self-adjoint.  Then  $W(\bA) \subseteq \IR^m$.
If $(\mu_1,\dots, \mu_m) \in \IR^m$ in {\rm (a)}, then
$A_j + \mu_j I$ is self-adjoint for $j = 1, \dots, m$, and
$W_k(A_1 + \mu_1 I, \dots, A_m + \mu_m I) \subseteq \IR^m$.
If $T$ is real in (b), then $B_1, \dots, B_m$ are self-adjoint,
and $W_k(B_1, \dots, B_p) \subseteq \IR^p$.
In (c) the operators $C_1, \dots, C_q$ can be chosen to be self-adjoint. The vector  $(s_1,\dots,s_m)$
and the matrix $R$ can be chosen to be real.

\item
One can determine the minimum value $q$ in Proposition \ref{affine} (c)
in terms of the affine dimension of $W(\bA)$, which is defined as the
minimum dimension of a real linear subspace $\bV$ of $\IC^m \equiv \IR^{2m}$
such that for a fixed row vectors $\tilde \bv \in \IC^m$,
$$W_k(A_1, \dots, A_m) =  W_k(H_1,G_1, \dots, H_m,G_m)
\subseteq \bV + \tilde \bv = \{\bv + \tilde \bv: \bv \in \bV\}.$$
\end{itemize}

In \cite{LP0} (4.5), it was shown that for $A \in\BH$, if
$\cH$ is an orthogonal sum of the  subspaces $\cH_1$ and $\cH_2$
and $A = B \oplus C \in \cB(\cH_1) \oplus \cB(\cH_2)$,
then
\begin{equation}\label{WABC}W_k(A) = \conv \cup\{W_{k_1}(B) + W_{k_2}(C):
0 \le k_j \le \dim \cH_j \hbox{ for } j = 1,2, \hbox{ and }
k_1 + k_2  = k \}.\end{equation}

We have the following generalization. The proof for general operators can
be  adapted from that of the finite dimensional case in
\cite[Theorem 4.4]{LPW}.

\begin{proposition} \label{direct-sum}
Let
$\bA = (A_1, \dots, A_m) \in \BH^m$
with $A_j = A_{j1} \oplus \cdots \oplus A_{jr} \in \cB(\cH_1)\oplus \cdots
\oplus \cB(\cH_r)$, where $\cH$ is an orthogonal sum of the subspaces
$\cH_1, \dots, \cH_r$. Then
$$W_k(\bA) \subseteq \conv \cW \subseteq \conv W_k(\bA),$$
where
$$\cW = \cup\{ W_{k_1}(\bA_1) + \cdots + W_{k_r}(\bA_r):
k_1 + \cdots + k_r = k, 0 \le k_j \le \dim \cH_j \hbox{ for }
j= 1, \dots, r\}$$
with $\bA_\ell = (A_{1\ell}, \dots, A_{m\ell})$ for $\ell = 1, \dots, r$, and as a convention,
$W_0(\bB) = \{(0, \dots, 0)\}$ for all $\bB\in \BH^m$.
In particular, the three sets are equal if $W_k(\bA)$ is convex.
\end{proposition}

It is known and not hard to see that for $A\in \BH$,
\begin{equation}\label{Wk+1}\frac{1}{k+1}W_{k+1}(A) \subseteq \frac{1}{k}W_k(A).
\end{equation}
We have the following extension in the absence of convexity.

\begin{proposition} \label{W(A)/k} Let $\bA \in \BH^m$. Then
$$\frac{1}{k+1}W_{k+1}(\bA) \subseteq \frac{1}{k}\conv W_k(\bA)
\quad \hbox{ so that } \quad
\frac{1}{k+1}\conv W_{k+1}(\bA) \subseteq \frac{1}{k}\conv W_k(\bA).$$
\end{proposition}

\section{Convexity}
\setcounter{equation}{0}

To study the convexity of $W_k(\bA)$, we can focus on $\bA \in \SH^m$ by
Proposition \ref{H+iG}.
If $\dim \cH = 2$, one can apply  Proposition \ref{affine} (c) to  show that
$W(A_1, \dots,  A_m)$ is determined by
$W(\tilde C_1, \dots, \tilde C_q)$ with $q \le 3$
and
$\tilde C_1, \dots, \tilde C_q \in \{X, Y, Z\}$ with
\begin{equation}
\label{pauli}
X = \begin{pmatrix} 0 & 1 \cr 1 & 0 \cr \end{pmatrix}, \quad
Y = \begin{pmatrix} 0 & -i \cr i & 0 \cr \end{pmatrix}, \quad
Z = \begin{pmatrix} 1 & 0 \cr 0 & -1 \cr \end{pmatrix}.
\end{equation}
The matrices $X, Y, Z$ are known as the Pauli's matrices in  physics literature.
It is known and easy to check that
$$W(X) = W(Y) = W(Z) = [-1,1],$$
$$W(X,Y) = W(Y,Z) = W(Z,X) = \{(a,b): a, b \in \IR, a^2+b^2 \le 1\}$$
is the unit disk, and
\begin{equation}
\label{sphere}W(X,Y,Z) = \{(a,b,c): a,b,c \in \IR, a^2+b^2 + c^2 = 1\}
\end{equation}
is the unit sphere.
The above observations were used to study $W_k(\bA)$
for $\bA \in M_n^m$ in \cite{LP1}, where $M_n^m$ is the set of all
$m$-tuples of $n\times n$ matrices. One can extend the results therein
to infinite dimensional operators.

\begin{proposition} \label{P3.1} Let $\bA = (A_1, \dots, A_m) \in \SH^m$.
\begin{itemize}
\item[{\rm (a)}] Suppose $\dim \cH = 2$. Then $W_k(\bA)$ is convex if and only if
the span of $\{I, A_1, \dots, A_m\}$ has dimension at most $3$.
\item[{\rm (b)}] Suppose $\dim \cH\ge 3$. If
the span of $\{I, A_1, \dots, A_m\}$ has dimension at most 4, then
$W_k(\bA)$ is convex. In particular, $W_k(A_1,A_2,A_3)$ is always convex.
\item[{\rm (c)}] Suppose $\dim \cH \ge 3$. If
the span of $\{I, A_1, \dots, A_m\}$ has dimension at least
4, then there is
$A_0 \in \BH$
such that $W_k(A_0, A_1,\dots, A_m)$ is not convex.
\end{itemize}
\end{proposition}

\it Proof. \rm
(a) In such a case, $W_k(\bA)$ is an affine transform of $W_k(C_1, C_2)$ for two $C_1, C_2\in \SH$
by Proposition \ref{affine} (c), which is always convex.

(b) Suppose $\dim \cH \ge 3$. If
the span of $\{I, A_1, \dots, A_m\}$ has dimension at most 4, then
as in (a),
$W_k(\bA)$ is an affine transform of $W_k(C_1,\dots, C_q)$ with
$q \le 3$. By the result in \cite{AP}
(see also \cite{AT}), $W_k(C_1, \dots, C_q)$ is convex,
and so is $W(\bA)$.

(c) Suppose $\span \{I, A_1, \dots, A_m\}$ has dimension 4 or larger.
By \cite[Theorem 4.1]{LP1},
there is an isometry $X: \IC^2 \rightarrow \cH$ such that
$\span \{I_2, X^*A_1X, \dots, X^*A_mX\}$ has dimension 4.
By renumbering the subscripts, we may assume that
$\{I_2, X^*A_1X, \dots, X^*A_3X\}$ is linear independent.
Let $P\in \BH$ be a rank $k-1$ orthogonal projection such that
$A_0 = 2P + XX^*\in \BH$ is a rank $k+1$ orthogonal projection.
We claim that $W_k(A_0,A_1,A_2,A_3)$ is not convex.
To see this, consider the set
$$\cS = W_k(A_0, \dots, A_3) \cap \{ (2k-1, a_1, a_2, a_3): a_1, a_2, a_3\in \IR\}.$$
If $(2k - 1, a_1, a_2, a_3)\in \cS$, then there is an orthogonal set
$\{x_1, \dots, x_k\}$ such that
$\sum_{j=1}^k   \la A_0x_j, x_j\ra = 2k-1$.
Suppose the compression of $A_0$ onto $\span \{x_1, \dots, x_k\}$
is $B_0 \in M_k$, which has eigenvalues $\beta_1 \ge \cdots \ge \beta_k$.
Since $A_0$ is unitarily similar to $2I_{k-1} \oplus I_2 \oplus 0$.
By interlacing inequalities, we see that
$2 \ge \beta_j$ for $j = 1, \dots, k-1$ and $1\ge \beta_k$.
Now, $\sum_{j=1}^k \beta_j = 2k-1$. So,
$\span \{x_1, \dots, x_k\}$ must contain all the eigenvectors of
$A_0$ corresponding to 2, and an eigenvector corresponding to the eigenvalue 1.
Note that if $\{x_1, \dots, x_k\}$ and $\{y_1, \dots, y_k\}$
are orthonormal sets with the same linear span, then
$\sum_{j=1}^k \la A x_j, x_j\ra =
\sum_{j=1}^k \la A y_j, y_j\ra$, we may simply assume that
$x_1, \dots, x_{k-1}$ are eigenvectors of $A_0$ corresponding to the eigenvalue 2,
and $x_k$ lies in the range space of $XX^*$.
If $(b_0, b_1, b_2, b_3) =
(\tr(P A_0), \tr(P A_1), \tr(P A_2), \tr(P A_3))$, then
$$(a_0, a_1, a_2, a_3) = (b_0, \dots, b_3)
+ (\la A_0x_k,x_k\ra, \dots, \la A_3x_k,x_k\ra).$$
Thus,
$$\cS = (b_0, \dots, b_3) + W(I_2, X^*A_1X, X^*A_2X, X^*A_3X),$$
where $W(X^*IX, X^*A_1X, X^*A_2X, X^*A_3X)
= W(I_2, X^*A_1X, X^*A_2X, X^*A_3X)$ is not convex by (a).
Hence, $\cS$ is not convex, and neither is $W_k(A_0, \dots, A_3)$.
\qed

\medskip
The next result extends Theorems 2.2 and 2.4 in \cite{LP1}. We use 0 to
denote the zero operator on some Hilbert space when confusion is unlikely.

\begin{proposition} \label{cS}
Let $\bA  = (A_1, \dots, A_m)\in \BH^m$. Suppose

\medskip
{\rm (1)} $\cS\subseteq \BH$ consists of operators of the form
 $aI + \begin{pmatrix} 0_k & R_1 \cr R_2 & 0\cr
\end{pmatrix}$, or

\medskip
{\rm (2)} $\dim \cH \ge 2k$ and
$\cS\subseteq \BH$
consists of operators of the form $R \oplus \alpha I$ with $R \in M_k,\ \alpha\in \IC$.

\medskip\noindent
If there is a unitary $U \in \BH$ such that
$U^*A_jU \in \cS$ for $j = 1, \dots, m$, then $W_k(\bA)$ is convex.
\end{proposition}

\it Proof. \rm
We may consider
$W_k(H_1,G_1, \dots, H_m,G_m)$ with
$A= H_j + iG_j$ with $H_j, G_j \in \SH$.
So, we may assume that $A_1, \dots, A_m \in \SH$.

If $\cS$ satisfies (1), we may assume that
$A_j =  a_jI + \begin{pmatrix} 0_k & R_j \cr R_j^* & 0_{n-k}\cr\end{pmatrix}$
for $j = 1, \dots, m$ so that $\bmu = (\mu_1, \dots, \mu_m) \in W_k(\bA)$
if and only if  there is a rank $k$ orthogonal projection
\begin{equation} \label{cQ}
Q = \begin{pmatrix} Q_{11} & Q_{12} \cr Q_{21} & Q_{22} \cr\end{pmatrix}
\in \BH \qquad \hbox{ with } Q_{11} \in M_k
\end{equation}
satisfying
$\mu_j = \tr (A_j Q) = ka_j + 2\re(\tr(R_j Q_{21}))$ for $j = 1, \dots, m.$
So,
$W_k(\bA) = k(a_1, \dots, a_m) + T(\cQ),$
where
$T$ is the real linear map
$X \mapsto 2(\re(\tr(R_1X)), \dots, \re(\tr(R_mX)))$
and $\cQ$ is the set of the $Q_{21}$ appeared as an $(2,1)$ block of a
rank $k$ orthogonal projection $Q$ in the form (\ref{cQ}).
Note that $Q_{21} \in \cQ$ if and only if $\|Q_{21}\| \le 1$.
This follows from the fact that
$Q_{21} \in \cQ$ if and only if there is a unitary $U = U_1 \oplus U_2\in \BH$
with $U_1 \in M_k$
and a nonnegative diagonal matrix $D = \diag(d_1, \dots, d_k)$
with $1 \ge d_1 \ge \cdots\ge d_k \ge 0$
such that
$$U^*QU =\begin{pmatrix} D & \sqrt{D-D^2} & 0 \cr
\sqrt{D-D^2} & I-D & 0 \cr
0 & 0 & 0 \cr\end{pmatrix}.$$
Thus, $T(\cQ)$ is convex.

The proof for  $\cS$ satisfying (2) is similar.
Assume $A_j = a_j I + (P_j \oplus 0)$ for $j = 1, \dots, m$.
Then
$$W_k(\bA) = k(a_1, \dots, a_m)  + \tilde T(\tilde \cQ),$$
where $\tilde T$ is the real linear map $X \mapsto
(\tr(P_1X), \dots, \tr(P_mX))$ and
$\tilde \cQ$ is the set of $Q_{11}$ appeared as the $(1,1)$ block of  a
rank $k$ orthogonal projection of the form $(\ref{cQ})$ so that
$\tilde\cQ$ is the convex set of positive semidefinite contractions in $M_k$. \qed

For any nonnegative integer $k$, we  present an example $\bA \in \BH^4$ such that
$W_{k+1}(\bA)$ is not convex but $W_r(\bA)$ is convex for all other positive integers
$r$ not larger than $\dim \cH$.
In particular, it shows that there are no implications between the
convexity of $W_k(\bA)$ and $W_{k+1}(\bA)$ for  $\bA \in \BH^m$ for $m \ge 4$.

\begin{example}  \label{3.3}
Let $k$ be a nonnegative integer.
Suppose $\dim \cH = N \ge 2(k+1)$, where $N$ could be infinity.
Let $\bA = (A_1, A_2, A_3, A_4)$ with $A_1 = X\oplus 0_{N-2}$,
$A_2 = Y\oplus 0_{N-2}$, $A_3
= Z \oplus 0_{N-2}$, and $A_4 = 0_2 \oplus I_{k} \oplus -I_{N-k-2}$, where
$X, Y, Z$ are the Pauli's matrices defined in (\ref{pauli}).
Then $W_{k+1}(\bA)$ is not convex.

{\rm (1)} If $r \le k$, then
$W_r(\bA)$ is convex and equals
$$\{(a,b,c,d) \in \IR^4: \xi = \sqrt{a^2+b^2 + c^2} \le 1,
\xi + |d| \le r\}.$$

{\rm (2)} If $k+1 < r$ and $2r$ is not larger than $N$, then
$W_r(\bA)$ is convex and equals
$$\{(a,b,c,d) \in \IR^4:
 \xi = \sqrt{a^2+b^2 + c^2} \le 1,  \hbox{ and }
-r \le d - \xi \le  2k - r + 1\}.$$

{\rm (3)} If $N$ is finite and $r \ge N/2$, then
$W_r(A) = (0,0,0, 2k+2-N) - W_{N-r}(\bA)$
so that $W_r(\bA)$ is convex except for $r = N-k-1$.
\end{example}

\it Proof. \rm Let $\{e_1, \dots, e_{2k+2}\} \cup \cB_0$ be an
orthonormal basis for $\cH$ so that $A_1,\dots, A_4$ have operator
matrices of the said form.

First, we show that $W_{k+1}(\bA)$ is not convex.
Consider the set
$W_{k+1}(\bA) \cap \{(a,b,c,k): a, b, c  \in \IR\}.$
If $(a,b,c,k) = \sum_{j=1}^{k+1} \la \bA x_j, x_j\ra$, then
$k = \sum_{j=1}^{k+1} \la A_4 x_j,x_j\ra$.
We may assume that
$(x_1, \dots, x_{k}) = (e_3, \dots, e_{k+2})$ and
$x_{k+1}  \in \span\{e_1, e_2\}$.
Then
 $$W_{k+1}(\bA) \cap \{(a,b,c,k): a, b, c  \in \IR\}
= \{(a,b,c,k): a, b, c \in \IR, a^2 + b^2 + c^2 = 1\}.$$
which is not convex. So, $W_{k+1}(\bA)$ is not convex.

(1) We show that $W_r(\bA)$ is convex if $r \le k$.
We prove the case when $r = k$, and the other cases are similar.
By Proposition \ref{direct-sum},  it suffices to show that
$W_k(\bA) = \conv \cW$,
where
$$\cW = \cup_{0\le \ell \le k} (W_\ell(\bB) + W_{k-\ell}(\bC)) \ \ \hbox{
with } \ \bB = (X,Y,Z,0_2) \hbox{  and } \ \bC = (0, 0, 0, I_k \oplus  -I_{N-k-2} ).$$
If $k = 0$, then $W_0(\bA) = \{(0,0,0,0)\}$ is convex.
Suppose $k = 1$. Then
$$\cW = W(X,Y,Z,0_2) \cup W(0,0,0,[1]\oplus -I_{N-3})$$
and
\begin{eqnarray*}
\conv\cW
&=&
\{ t(\xi_1,\xi_2,\xi_3,0) + (1-t)(0,0,0,\xi_4)\in\IR^4:
\xi_1^2 + \xi_2^2 + \xi_3^2 = 1, |\xi_4|\le 1\}
\\
&=& \{(a,b,c,d)\in\IR^4: \sqrt{a^2+ b^2 + c^2} + |d| \le 1 \}.
\end{eqnarray*}
\noindent
Suppose $k > 1$. Then

\medskip
$W_2(\bB)+ W_{k-2}(\bC) = \{(0,0,0,d_3): d_3 \in [2-k,k-2]\}$ is a subset of

\medskip
$W_0(\bB) + W_k(\bC) = \{(0,0,0,d_1): d_1 \in [-k,k]\}$, and

\medskip
$W_1(\bB) + W_{k-1}(\bC)
= \{(\xi_1, \xi_2, \xi_3, d_2)  \in \IR^4:
\xi_1^2 + \xi_2^2 + \xi_3^2 = 1,
d_2 \in [1-k, k-1]\}$.

\noindent
As a result, $\conv \cW$
consists of elements $\bmu =(1-t)(0,0,0,d_1)+ t(\xi_1,\xi_2,\xi_3, d_2)$
with
$$t \in [0,1], \ d_1 \in [-k,k], \ d_2 \in [1-k,k-1], \
(\xi_1,\xi_2,\xi_3) \in \IR^3 \hbox{ is  a unit vector},$$
i.e.,  $\bmu = (t\xi_1,t\xi_2,t\xi_3,d_t)$ such that
$(\xi_1, \xi_2, \xi_3) \in \IR^3$ is a unit vector and
$$d_t \in [(1-t)(-k) + t(1-k), (1-t)k + t(k-1)] = [-k+t, k-t].$$
So,
$$ \conv \cW  = \{(a,b,c,d)\in \IR^4:  a^2 + b^2 + c^2 \le 1,
\sqrt{a^2 + b^2 + c^2} + |d| \le k \}.$$

We will show that $\conv \cW \subseteq W_k(\bA)$ in the following.
For any $(a,b,c,d) \in \conv\cW$ with $\sqrt{a^2 + b^2 + c^2} = t$, let
$v \in \span\{e_1, e_2\}$ be a unit vector such that
$t(v^*A_1v, v^*A_2v, v^*A_3v) = (a,b,c)$.
For $s \in [0, \pi/2]$, let $Q_s$ be the rank $k$ orthogonal projections
with range space spanned by

\medskip\centerline{
$\sqrt t v + \sqrt{1-t} (\cos s \ e_3 + \sin s \ e_{k+3})$
\quad and \quad $\cos s \ e_j + \sin s \ e_{k+j}$ \  for
$j = 4, \dots, k+2$.}

\medskip\noindent
Here if $k = 1$, we only have the first vector in the above list.
Then
$$(\tr(A_1Q_s), \dots, \tr(A_4Q_s))
= (a,b,c, d_s) \in W_k(\bA),$$
with
$d_s = ((1-t) + (k-1))\cos^2s-((1-t) + (k-1))\sin^2 s=
 (k-t)\cos 2s
  \in W_k(\bA)$
covering all values $d$ with
$\sqrt{a^2+b^2 + c^2} + |d| \le k$.
Thus, $W_k(\bA) = \conv \cW$ is convex.

(2) Suppose $r > k+1$ and and $2r$ is not larger than $N$.
By Proposition  \ref{cS} (2) with $A_4$ replaced by
$A_4+I_N$, and $k$ replaced by $r \ge k+2$, we see that
$W_r(\bA)$ is convex. By Proposition \ref{direct-sum},
$W_r(\bA) = \conv \cW$ with
$$\cW = \cup_{0\le \ell \le k} (W_\ell(\bB) + W_{k-\ell}(\bC)) \ \ \hbox{
with } \ \bB = (X,Y,Z,0_2) \hbox{  and } \ \bC = (0, 0, 0, I_k \oplus  -I_{N-k-2})).$$
If $k = 1$, then
$$\cW = W(X,Y,Z,0) \cup W(0,0,0,[1]\oplus -I_{N-3})$$
and
\begin{eqnarray*}
\conv\cW
&=&
\{ t(\xi_1,\xi_2,\xi_3,0) + (1-t)(0,0,0,\xi_4)\in\IR^4:
\xi_1^2 + \xi_2^2 + \xi_3^2 = 1, |\xi_4|\le 1\}
\\
&=& \{(a,b,c,d)\in\IR^4: \sqrt{a^2+ b^2 + c^2} + |d| \le 1 \}.
\end{eqnarray*}
Suppose $k > 1$. Then

\medskip
$W_2(\bB)\oplus W_{r-2}(\bC) = \{(0,0,0,d_3): d_3 \in
[-r+2,k - (r-k-2)]\}$,

\medskip
$W_0(\bB) + W_r(\bC) = \{(0,0,0,d_1): d_1 \in [-r,k-(r-k)]\}$, and

\medskip
$W_1(\bB) + W_{r-1}(\bC)
= \{(\xi_1, \xi_2, \xi_3, d_2): \xi_1, \xi_2, \xi_3, d_2 \in \IR,$

\medskip
$\hskip 2in
\xi_1^2 + \xi_2^2 + \xi_3^2 = 1,
d_2 \in [-r+1, k - (r-k-1)]\}$.

\medskip\noindent
One readily gets the description of
$\conv \cW = W_r(\bA)$.

\medskip
(3) The assertion follows from the fact that $W_\ell(\bA) = (\tr A_1, \dots, \tr A_4)-
W_{N-\ell}(\bA)$ if $N$ is finite.
\qed

\medskip
Note that when $k = 1$, the matrices
$A_1, \dots, A_4\in M_4$
in Example \ref{3.3} have the lowest dimension
satisfying $W(\bA)$ is convex and $W_2(\bA)$ is not.
For $\bA=(A_1,\dots, A_4) \in M_3^4$, we have
$W(\bA) = (\tr A_1, \dots, \tr A_4) - W_2(\bA)$ so that
$W(\bA)$ is convex if and only if $W_2(\bA)$ is convex.

\section{Closure of $W_k(\bA)$ and $\conv W_k(\bA)$}
\setcounter{equation}{0}

In the proof of Theorem \ref{P3.1},
the set $\cS$ is  the intersection of
$W(A_0, \dots, A_3)$ and a support plane $\{(2k-1, a_1, a_2, a_3):
a_1, a_2, a_3 \in \IR\}$ of $\conv W_k(A_0, \dots, A_3)$.
In the verification of Example \ref{3.3}, the set
$W(A_1, A_2, A_3, A_4) \cap \{(a,b,c,k): a,b,c \in \IR\}$ is
the intersection of $W(A_0, \dots, A_3)$ and a support plane of
$\conv W_k(A_1, \dots, A_4)$.
In general, we can use the spectrum $\sigma(X)$ of a self-adjoint
operator $X$ to help describe the support planes of
$\conv\cl (W_k(\bA)) = \cl(\conv W_k(\bA))$, and then express
$\cl (\conv W_k(\bA))$ as the intersection of the half spaces. This
result is useful for the study of the geometric feature of $W_k(\bA)$ and
$\conv W_k(\bA)$.

To present the construction, we need some notations and a lemma.
Denote by $\sigma(X)$ the spectrum of $X \in \BH$,
For any $A \in \SH$ and $\ell \in \IN$, define
$$\lambda_\ell(A) = \sup \{\|A-X\|: X\in \SH \hbox{ has rank }
\ell-1\}.$$

\begin{lemma} \label{s-number}
Let $A \in \SH$ where $\dim \cH$ is infinite.
Then $A$ is unitarily similar to $A_1 \oplus A_2$, for a diagonal
operator $A_1$ and one of $A_1$ or $A_2$ may be void, with
one of the following holds:
\begin{itemize}
\item[{\rm (a)}] For some finite $\ell \ge 1$,
$A_1 = \diag(\lambda_1(A),
\dots, \lambda_\ell(A))$,
$\mu=\max\sigma(A_2) \le \lambda_\ell(A_1)$
and $\mu$ is not an eigenvalue of $A_2$; hence
$\lambda_p(A) = \mu$ for all $p > \ell$.
\item[{\rm (b)}] The operator $A_1 = \diag(\lambda_1(A), \lambda_2(A), \dots )$
has infinite dimensional range and
$$\mu=\max\sigma(A_2) \le \inf \{\lambda_\ell(A_1) :\ell\ge 1\}.$$
\end{itemize}
\end{lemma}

One readily sees that for $A\in \SH$,
$$\cl(W_k(A)) = \big[-\sum_{j=1}^k\lambda_j(-A), \sum_{j=1}^k \lambda_j(A)\big].$$
For $W_k(\bA)$, we have the following.

\begin{proposition} \label{halfspaces}
Let  $W_k(\bA)$ with $\bA = (A_1, \dots, A_m) \in \SH^m$.
Suppose  $\bu = (u_1, \dots, u_m)$ is
a real unit vector and $A_\bu = \sum_{j=1}^m u_j A_j$.
Then
\begin{equation}\label{halfspace}
P_\bu(k,\bA) = \left\{(\zeta_1,\dots, \zeta_m):
\sum_{j=1}^m u_j \zeta_j  \le \sum_{j=1}^k \lambda_j(A_\bu)\right\}
\end{equation}
is a half space containing $W_k(\bA)$ and
the boundary of $P_\bu(k,\bA)$
\begin{equation}\label{supportplane}
\partial P_\bu(k,\bA)  = \left\{(\zeta_1,\dots, \zeta_m):
\sum_{j=1}^m u_j \zeta_j  = \sum_{j=1}^k \lambda_j(A_\bu)\right\},
\end{equation}
is a support plane of $W_k(\bA)$.
Consequently,
$$\cl(\conv W_k(\bA)) = \cap\{ P_\bu(k,\bA): \bu \in \IR^m \hbox{ a unit vector} \}.$$
\end{proposition}

In \cite{CLP0},
we showed that
$\cl(W_k(A))= \conv \cS$, where
$$\cS   = \bigcup_{0 \le \ell \le k} (W_\ell(A) + (k-\ell)W_{ess}(A)).$$
Equivalently, the extreme points of $\cl(W_k(A))$ lie in $\cS$.
As a result,
$W_{k}(A)$ is closed if and only if $\cS \subseteq W_k(A)$.
Here, recall that we use the convention that $W_0(A) = \{0\}$.
This extends the results of Lancaster \cite{L} for a single operator $A \in \BH$.
Moreover, we showed in \cite{CLP0} that if $W_{k+1}(A)$ is closed, then so is $W_k(A)$.
The proofs therein depend on the convexity of $W(A)$.

In the following, we
show that the results can be extended to $\conv W_k(\bA)$ for $\bA \in \BH^m$,
and the results will be useful for studying the polyhedral properties of $W_k(\bA)$
and $\conv W_k(\bA)$ in later sections.
We need the concept of the joint essential numerical of $\bA \in \BH^m$ defined
in (\ref{Wess}). There are several
equivalent characterizations of $W_\ess(\bA)$, see \cite{LP2}.
For our purpose, we will use the fact that
a point $\bmu= (\mu_1, \dots, \mu_m)\in \mathbb C^m$ belongs to $W_\ess(\bA)$
if and only if there is a weakly null sequence of unit vectors (or,
a sequence of  orthonormal vectors) $\{v_k\}$ in $\cH$ such that
$\langle A_jv_k, v_k \rangle \to \mu_j$ for $1\le j\le m$.
We will continue to use the convention that
$W_0(\bA) = \{(0,\dots, 0)\}$
for any $\bA \in \BH^m$.

\begin{theorem} \label{4.3}
Let $\bA \in \SH^m$. Consider the Hilbert space $\cH_1$, which may be inseparable,
with an orthonormal basis
$\{e_\bmu: \bmu= (\mu_1, \dots, \mu_m) \in W_{\rm ess}(\bA)\}$.
Let $(D_1, \dots, D_m)$ be the $m$-tuple of diagonal operators
acting on $\cH_1$ by $D_j e_\bmu = \mu_j e_\bmu$, and
$\cH_0 = \cH_1\oplus \cdots\oplus \cH_1$ be the orthogonal sum
of $k$ copies of $\cH_1$. Suppose
$\tilde \bA = (\tilde A_1, \dots, \tilde A_m)\in\( \cS(\cH)  \oplus \cS(\cH_0)\)^m$
with $\tilde A_j =  A_j \oplus  (I_k\otimes D_j)$ for $j = 1, \dots, m$,

\begin{itemize}
\item For $\ell = 1, \dots, k$,
$W_\ell (I_k\otimes D_1, \dots, I_k\otimes D_m) = \ell W_{\ess}(\bA)$.
\item The equality $\cl W_k(\bA) = W_k(\tilde \bA)$ holds.
\item
The set $W_k(\bA)$ is closed if and only if $W_k(\bA) = W_k(\tilde \bA)$,
equivalently, every element in $W_k(\bA)$ has the form
$\(\tr(A_1P),\dots,\tr(A_mP)\) + (k-(\tr P))\bxi$ for a positive semidefinite
contraction $P$ with rank at most $k$
and $\bxi \in  W_{\ess}(\bA)$.
\end{itemize}
\end{theorem}

\it Proof. \rm
By \cite[Theorem 3.1]{LP2}, $W_\ess(\bA)$ is convex and closed.
Thus, $W(D_1, \dots, D_m) = W_\ess(\bA)$.
Let $\bD = (I_k \otimes D_1, \dots, I_k \otimes D_m)$.
Then $W(\bD) = W(D_1, \dots, D_m)$.
Suppose $\ell \in \{2, \dots, k\}$.
If $\bmu \in W(D_1, \dots, D_m)$, then
$\ell \bmu \in \ell W(D_1, \dots, D_m)\subseteq W_\ell(\bD)$.
If $\{x_1, \dots, x_\ell\}$
is an orthonormal set, and $(a_1, \dots, a_m) \in
W_\ell(\bD)$
with $a_j = \sum_{i=1}^\ell \la (I_k\otimes D_j) x_i, x_i\ra$ for $j = 1, \dots, m$.
Since
$$(\la (I_k\otimes D_1) x_i,x_i\ra, \dots, \la (I_k\otimes D_m) x_i,x_i\ra) \in W(\bD)
= W(D_1, \dots, D_m),$$
it follows that
$$(a_1, \dots, a_m) \in \{\bxi_1 + \cdots + \bxi_\ell:
\bxi_i \in W(D_1, \dots, D_m), i = 1, \dots, \ell \}
= \ell W(D_1, \dots, D_m),$$
where the set equality holds because $W(D_1,\dots, D_m)$ is closed and convex.

For the second assertion, we first show that $\cl(W_k(\bA)) \subseteq W_k(\tilde \bA)$.
Let $\bmu\in \cl(W_k(\bA))$.  There are orthonormal sets of vectors
$\{v^{(n)}_1, \ldots, v^{(n)}_k\}$ in $\mathcal H$ such that
$\sum_{j = 1}^k \langle \bA v^{(n)}_j, v^{(n)}_j\rangle \to\bmu$.
As the closed unit ball of $\cH$ is weakly sequentially compact, by
passing to subsequences, we may assume that
for each $j$, $v^{(n)}_j\to v_j$ weakly and
$\langle \bA v^{(n)}_j, v^{(n)}_j\rangle\to \bmu_j$,
for $v_j$ in the closed unit ball of $\mathcal H$ and
$\bmu_j\in \cl(W(\bA))$.
There are three possibilities (see the proof of
\cite[Theorem 2.1]{C} for detail),
\begin{enumerate}
\item[(i)] $v_j = 0$ and $\bmu_j\in W_\ess(\bA)$,

\item[(ii)] $\|v_j\| = 1$, $v_j^{(n)}\to v_j$ strongly and
$\bmu_j = \langle \bA v_j, v_j\rangle \in W(\bA)$,

\item[(iii)] $0 < \|v_j\| < 1$ and
$\bmu_j =\|v_j\|^2\big\langle \bA\frac {v_j}{\|v_j\|},
\frac {v_j}{\|v_j\|}\big\rangle + (1 - \|v_j\|^2)\bxi_j$
for some $\bxi_j\in W_\ess(\bA)$ so that $\bmu_j$ is a convex
combination of points in $W(\bA)$ and $W_\ess(\bA)$.
\end{enumerate}

Taking any $\bxi_j\in W_\ess(\bA)$ in (ii), we can always write
$\bmu_j = \la \bA v_j, v_j\ra + (1 - \|v_j\|^2)\bxi_j$.
As in \cite{LP0}, consider the positive semidefinite operator
$H = \sum_{j = 1}^k \langle\,\cdot\,, v_j\rangle v_j$.
Let $d_1\ge d_2\ge\cdots \ge d_k\ge 0$ be the $k$ largest eigenvalues
of $H$, and $\{u_1, \ldots, u_k\}$ an orthonormal set of corresponding
eigenvectors. Then for each $j$,
$$d_j = \langle Hu_j, u_j\rangle
= \sum_{i = 1}^k |\langle u_j, v_i\rangle|^2
= \lim_{n\to\infty} \sum_{i = 1}^k |\langle u_j, v_j^{(n)}\rangle|^2
\le \|u_j\|^2 = 1,$$
and $d_1 + \cdots + d_k = \mbox{tr} H
= \sum_{j = 1}^k \|v_j\|^2$.
We also have
$$\sum_{j = 1}^k (1 - \|v_j\|^2)\bxi_j = (k - \sum_{j=1}^k d_j)
\left(\frac {1 - \|v_1\|^2}{k - d_1 - \cdots - d_k}\bxi_1 + \cdots
+ \frac {1 - \|v_k\|^2}{k - d_1 - \cdots - d_k}\bxi_k\right).$$
Since $W_\ess(\bA)$ is convex, see for example \cite{LP2},
$\sum_{j = 1}^k (1 - \|v_j\|^2)\bxi_j = (k - \sum_{j=1}^k d_j)\bxi$
for some $\bxi\in W_\ess(A)$.
Thus, $\bmu = (\tr(A_1H), \dots, \tr(A_mH)) + (k-\tr H) \bxi$.
Applying a suitable unitary similarity transformation to $H, A_1, \dots, A_m$, we may
assume that $H = T \oplus 0$ with $T = \diag(d_1,\dots, d_k)$.
Applying a suitable unitary similarity transform to $\bD$, we may
assume that  the first $k$
diagonal entries in $\bD$  are equal  to $\bxi$.
Let
 $P = \begin{pmatrix} P_{11} & P_{12} \cr P_{21} & P_{22} \cr\end{pmatrix}
 \in \cB(\cH \oplus \cH_0)$ be
 such that

\medskip\centerline{
 $P_{11} = T \oplus 0 \in \SH$,
 $P_{22} = (I_k-T) \oplus 0 \in \cS(\cH_0)$
and $P_{12} = \sqrt{T-T^2} \oplus 0 = P_{21}^*$.}

\medskip\noindent
Then
$P$ is a rank $k$ orthogonal projection, and
$\bmu = (\tr(\tilde A_1P), \dots, \tr(\tilde A_mP)) \in W_k(\tilde \bA).$

\medskip
To prove the reverse inclusion,
suppose $P$ is a rank $k$ orthogonal projection on $\cB( \cH \oplus \cH_0)$
and
$\bmu = (\tr(\tilde A_1P), \dots, \tr(\tilde A_mP)) \in W_k(\tilde \bA).$
Let $P = \begin{pmatrix} P_{11} & P_{12} \cr P_{21} & P_{22} \cr\end{pmatrix}$
such that $P_{11} \in \cB(\cH)$ and $P_{22} \in \cB(\cH_0)$.
We will show that $\bmu\in \cl W_k(\bA)$.

First, note that
there is a unitary $U =U_1 \oplus U_2\in \cB(\cH \oplus \cH_0)$ such that
$U_1^*P_{11}U_1 = T \oplus 0\in \cS(\cH)$ with $T = \diag(d_1, \dots, d_k)$ and
$U_2^*P_{22}U_2 = (I-H) \oplus 0 \in \cS(\cH_0)$.
We may assume that $U_1 = I_\cH,U_2 = I_{\cH_0}$ by replacing
$X$ by $U^*XU$ for all $X \in \{P, \tilde A_1, \dots \tilde A_m\}$,
where $\tilde A_j = A_j \oplus \tilde D_j$. Here note that
$\tilde D_j = U_2^*(I_k\otimes D_k)U_2$ is no longer of the form $I_k \otimes D_j$,
for $j = 1, \dots, m$, but
we still have  $W(\tilde D_1, \dots, \tilde D_n) = W(I_k \otimes D_1, \dots,
I_k\otimes D_m) = W_\ess(\bA)$.
Suppose $B_j$ is the leading $k\times k$ submatrix of $A_j$
and $C_j$ is the leading $k\times k$ submatrix of $\tilde D_j$.
Then $\mu_j = \tr (B_j T) + \tr (C_j(I_k-T))$ for $j = 1, \dots, m$.
For any $\varepsilon > 0$, we can find a unitary
$V = I_k \oplus V_0 \in \cB(\cH)$ such that for $j = 1, \dots m$,
$V^*A_jV =
\begin{pmatrix} B_j & F_j & 0 \cr
               F_j^* & \star &  \star\cr
                0 & \star & \star \cr\end{pmatrix}$,
where $F_j$ is an $k\times mk$ matrix because
the $(1,2)$ block to $A_j$ has rank at most $k$ for each $j$.
Since the diagonal entries of $(\tilde D_1, \dots, \tilde D_m)$ lies in
$W_{\rm ess}(\bA)$, there is a unitary $W = I_{(m+1)k} \oplus W_0$ such that
for $j = 1, \dots, m$,
$W^*V^*A_jVW =
\begin{pmatrix} B_j & F_j & 0 \cr
               F_j^* & \star &  \star\cr
                0 & \star & L_j \cr\end{pmatrix}$,
where the $\ell_1$-norm between the vector of first
$k$ diagonal entries of $L_j$ and those of $C_j$
is bounded by $\varepsilon$. This can be done because
the $m$-tuples of diagonal entries of $(C_1, \dots, C_m)$
lie in $W_{\rm ess}(\bA)$. Let
$\tilde P = \begin{pmatrix}
\tilde P_{11} & \tilde P_{12} \cr
\tilde P_{21} & \tilde P_{22} \cr \end{pmatrix}$
with
$\tilde P_{11} = T \oplus 0_{mk}$,
$\tilde P_{22} = (I_k-T) \oplus 0$,
and $\tilde P_{12} = \sqrt{T-T^2} \oplus 0 = \tilde P_{21}^*$.
Then
$$|\bmu - (\tr(A_1\tilde P), \dots, \tr(A_m \tilde P))|
= |(\tr(\tilde P_{22}(C_1\oplus 0-L_1)), \dots, \tr(\tilde P_{22}(C_m\oplus 0-L_m)))|
< \varepsilon.$$
So, the second assertion is valid.

For the third assertion, it is clear
that $W_k(\bA)$ is closed if and only if $W_k(\bA) = W_k(\tilde \bA)$ by the
second assertion. Evidently, every element in
$W_k(\tilde \bA)$ has the form
$(\tr(\tilde A_1 Q), \dots, \tr(\tilde A_m Q))$ for a rank $k$ orthogonal projection $Q$.
We may assume that
 $Q = \begin{pmatrix}X_1X_1^*&X_1X_2^*\cr X_2X_1^*&X_2X_2^*\end{pmatrix}   $ where $X_1:\IC^k\to \cH$ and $X_1:\IC^k\to \cH_0 $ satisfy $X_1^*X_1+X_2^*X_2=I_k$.  Then $P=X_1X_1^*$ is a positive semidefinite contraction  in  $\SH$ with rank at most $k$. There is a unitary operator $V \in \cB(\cH_0)$ such that $X_2X_2^* = VFV^*$ with $F = \diag(f_1, \dots, f_k) \oplus 0$
 such that $\sum_{j=1}^k f_j = k-\tr P$.
Since $W(D_1, \dots, D_m) = W(I_k\otimes D_1, \dots, I_k\otimes D_m)$,
we see that
$$
(\tr((I_k\otimes D_1)X_2X_2^*), \dots, \tr((I_k\otimes D_m)X_2X_2^*))
= \sum_{j=1}^k f_j \bxi_j$$
for some $\bxi_1, \dots, \bxi_k \in W_{\ess}(\bA)$.
Now, $W_{\ess}(\bA) = W(D_1, \dots, D_m)$ is convex, it follows that
$$\sum_{j=1}^k f_j \bxi_j = (\sum_{j=1}^k f_j) \bxi = (k-\tr P) \bxi$$
for some $\bxi\in W_{\ess}(\bA)$.
Thus,
 $(\tr(\tilde A_1 Q), \dots, \tr(\tilde A_m Q))$ equals
$$(\tr(  A_1 P), \dots, \tr(  A_m P))+
 (\tr((I_k\otimes D_1)X_2X_2^*), \dots, \tr((I_k\otimes D_m)X_2X_2^*)).$$
So, every element in $W(\bA)$ has the asserted form.
\qed

\begin{theorem} \label{closure}
Let $\bA = (A_1,\dots, A_m) \in \SH^m$.
Then
$$\cl(\conv W_k(\bA)) = \conv(\cl(W_k(\bA)))= \conv \cS,$$
where
$$\cS =
 \cup_{0 \le \ell \le k}
(W_\ell(\bA) + (k-\ell) W_{\ess}(\bA)).$$
Consequently, $\conv W_k(\bA)$ is closed if and only if
\begin{equation}\label{wk}W_\ell(\bA) + (k-\ell) W_{\ess}(\bA)
\subseteq \conv W_k(\bA)\quad \mbox{for}\quad \ell = 0, \dots, k-1.
\end{equation}
\end{theorem}

\it Proof. \rm The equality $\cl(\conv W_k(\bA)) =
\conv(\cl(W_k(\bA)))$ follows from the fact that $W_k(\bA)$
is bounded in $\mathbb R^m$ and Carath\'eodory's theorem.

For the second equality, we have $\conv(\cl(W_k(\bA)))= \conv W_k(\tilde \bA)$,
where $\tilde \bA$ is the operator defined in Theorem \ref{4.3}. Using the construction
there with $\bD = (D_1, \ldots, D_m)$ and  Proposition \ref{direct-sum}, we see that
$\conv W_k(\tilde \bA)$ is the convex hull of
$$\cup\{W_{r_1}(\bA) + W_{r_2}(\bD) + \cdots + W_{r_{k + 1}}(\bD): r_j\ge 0 \hbox{ for all $j$ and }
r_1 + \cdots + r_{k + 1} = k\}.$$
As $W(\bD) = W_{\ess}(\bA)$ is convex, it is not hard to see that
$W_{r_2}(\bD) + \cdots + W_{r_{k + 1}}(\bD) = (k - r_1)W_{\ess}(\bA)$.
The result follows.
\qed

If $\bA \in \SH^m$ with $m \le 3$,
 then $W_k(\bA)$ is convex. In such a case,
 $W_k(\bA)$ is closed if and only if
$W_\ell(\bA) + (k-\ell) W_{\ess}(\bA) \subseteq W_k(\bA)$
for $\ell = 0, \dots, k-1$ by  Theorem \ref{closure}.
Suppose $\bA \in \SH^m$ with $m > 3$.
If $W_k(\bA)$ is closed, then $\conv W_k(\bA)$ is closed by
Carath\'eodory's theorem.
Hence (\ref{wk}) holds.  However, $W_k(\bA)$ may not be closed even
if $\conv W_k(\bA)$ is closed, as shown in Example \ref{4.6},
which is a modification of \cite[Example 3.2]{LP2}.
In particular, Example \ref{4.6} shows that there is $\bA \in \KH^4$
such that $W(\bA)$ is not closed
even if $\conv W(\bA)$ is closed.
It is interesting to note that for a single compact operator $A \in \KH$,
$W(A)$ is closed if and only if $0 \in W(A)$ by  Theorem \ref{closure}
as $W_{\ess}(A) = \{0\}$.
For $W(\bA)$ with $\bA \in \KH^m$, $W(\bA)$ is closed if and only if
$W(\bA) = W(\bA \oplus (\b0, \dots, \b0))$ and
$\conv W(\bA)$ is closed if and only if $\{(0,\dots, 0)\} \in W(\bA)$.

\begin{example} \label{4.6}
Let $\bA = (A_1, \dots, A_4)$ with
$A_1=\diag(1,0,1/3,1/4,\dots), A_2 =\diag(1,0)\oplus 0,$
$$A_3=\begin{pmatrix} 0&1\cr 1&0\end{pmatrix} \oplus 0 \quad \hbox{ and }
\quad A_4=\begin{pmatrix} 0& -i\cr i&0\end{pmatrix} \oplus 0.$$
Then $W(\bA)$ is not closed but $\conv W(\bA)$ is closed.
\end{example}

\it Proof. \rm Since $W_{\ess}(\bA) = \{(0,0,0,0)\} \subseteq W(\bA)$,
$\conv W(\bA)$ is closed by Theorem \ref{closure}.
In \cite{LP2} that
it was shown that  $(1,1,0,0)/2 \in
 W(\bA \oplus (\b0, \dots, \b0)) = \cl(W(\bA))$, but
 $(1,1,0,0)/2 \notin W(\bA)$.
 Thus,  $W(\bA)$ is not closed.
\qed

One may also wonder whether $kW_\ess (\bA)\subseteq W_k(\bA)$ alone would imply
$W_k(\bA)$ or $\conv W_k(\bA)$ is closed. Example \ref{4.6} shows that, unlike
the single operator case, $W(\bA)$ may not be closed even if $W_\ess (\bA)\subseteq W(\bA)$.
This justifies our consideration of the convex hull $\conv W_k(\bA)$.
Also, let $A = 0_{k-2} \oplus \diag(1, -1, -1/2, -1/3, \dots)$ with $k \ge 2$.
Then $kW_{\rm ess}(A) = W_{\rm ess}(A) = \{0\}$ and
$0 = \tr (AP) \in W_k(A)$ if $P = I_k \oplus 0$.
However, $\sup W_k(A) = 1\notin W_k(A)$ so that $W_k(A) = \conv W_k(A)$ is not closed.

Note that there are no implication relations between the
convexity of $W_k(\bA)$ and $W_{k+1}(\bA)$ as shown in the last section.
It is somewhat surprising that we have the following.

\begin{theorem} \label{k-closed}
Let $\bA \in \SH^m$, and $k$ be a positive integer.
If $m \le 3$ and $W_{k+1}(\bA)$
is closed, then so is $W_k(\bA)$.
For $m \ge 4$, if $\conv W_{k+1}(\bA)$
is closed (in particular, if $W_{k+1}(\bA)$ is closed), then $\conv W_{k}(\bA)$ is closed.
\end{theorem}

The rest of this section is devoted to the proof of Theorem \ref{k-closed}.
Readers may skip the long technical proof and move to the next section directly.

\medskip
\it Proof of Theorem \ref{k-closed}. \rm We only need to prove the second assertion
for $m \ge 1$. The first assertion will then follow
as $W_\ell(\bA)$ is convex for any positive integer $\ell$
if $m \le 3$; see \cite{AT}.

To prove the second assertion, assume that $\conv W_{k+1}(\bA)$ is closed.
The set $\conv W_k(\bA)$ is closed if and only if every extreme point
of $\cl (\conv W_k(\bA))$ lies inside $\conv W_k(\bA)$.

Let $\bnu = (\nu_1, \dots, \nu_m)$ be an
extreme point of  $\cl (\conv W_{k}(\bA))$.
By Theorem 4.4, there is $r \in \{0, \dots, k\}$
and an orthonormal family $\{z_1, \dots, z_r\} \subseteq \cH$
such that
\begin{equation}\label{bxi}
\bnu = \bmu + (k-r) \bxi
\end{equation}
with
$$
\bmu = (\mu_1, \dots, \mu_m) = \sum_{j=1}^r \la\bA z_j, z_j\ra
\in W_r(\bA)
\ \hbox{and} \
\bxi = (\xi_1, \dots, \xi_m)  \in W_{\ess}(\bA).$$
We may assume that
$r$ is the {\bf maximum number} in $\{0, 1, \dots, k\}$
such that (\ref{bxi}) holds.

If $r = k$, then $\bnu \in W_k(\bA)$, and we are done.
So, assume that $r < k$. Then $\bnu \notin W_k(\bA)$.
We will derive a contradiction
by showing that $\bnu \in W_k(\bA)$.

Let $\bu = (u_1, \dots, u_m)\in \mathbb R^m$ be a unit vector that
such that $$\partial P_\bu(k,\bA)  = \big\{(\zeta_1,\dots, \zeta_m):
\sum_{j=1}^m u_j \zeta_j  = \sum_{j=1}^k \lambda_j(A_\bu)\big\}$$ is a
support plane of $\cl (\conv W_{k}(\bA))$ at $\bm \nu$.
See Proposition \ref{halfspaces}.
We can do the following.

\noindent
{\bf Reduction} \it
We may permute the entries of
$(u_1, \dots, u_m)$ and assume $u_m \ne 0$, and
permute the components of $(A_1, \dots, A_m)$  accordingly.
Then replace
$(A_m, \nu_m)$  by
$(\sum_{j=1}^m u_jA_j, \sum_{j=1}^m u_j\nu_j).$
After the replacement
${\bm \nu}$ remains to be an extreme point of the set
$\cl (\conv W_k(\bA))$.
Moreover, ${\bm \nu}
\in W_k(\bA)$ in the original problem
if and only if the same conclusion holds
for the modified ${\bm \nu}$ and $\bA$.
\rm

After the replacement,  we have
$(u_1, \dots, u_m) = (0,\dots, 0,1)$ and
$$\cl (\conv W_k(\bA)) \cap \partial P_\bu(k, \bA)\\
=  \left\{(\zeta_1, \dots, \zeta_m)\in \cl(\conv W_k(\bA)): \zeta_m =
 \sum_{j=1}^k  \lambda_j(A_m) \right\}.$$
It follows that
$$\la A_m z_j , z_j \ra \ge \xi_m \in W_{\ess}(A_m)
\qquad \hbox{ for } j = 1, \dots, r.$$
In particular, $\lambda_1(A_m), \dots, \lambda_r(A_m)$ are the
eigenvalues of the compression of $A_m$ onto ${\rm span}\{z_1, \dots, z_r\}$.
Applying a unitary transform, we may assume that
$$A_m z_j = \lambda_j(A_m) z_j, \qquad  j = 1, \dots, r,$$
and
$$\lambda_{r+1}(A_m) = \cdots =  \lambda_{k}(A_m) =\xi_m= \max W_{\ess}(A_m).$$

Let $r_1 = 0$ if $\lambda_1(A_m) = \xi_m$, and
$r_1 = \max \{ j: \lambda_j(A_m)  > \xi_m\}$ otherwise.
The vectors $z_1, \ldots, z_{r_1}$ are uniquely determined up to a unitary transform.
If $\{x_1,\dots, x_k\}\subseteq \cH$ is
an orthonormal family and $\sum_{j=1}^k \la \bA x_j, x_j \ra
\in W_k(\bA) \cap \partial P_\bu(k,A)$, we may assume that
$\la A_m x_1, x_1\ra \ge \cdots \ge \la A_m x_k, x_k\ra$, and
then further assume that $(x_1, \dots, x_{r_1}) = (z_1, \dots,
z_{r_1})$.
Similarly, we have
$W_{k+1}(\bA) \subseteq P_\bu(k+1,\bA)$
with
$$P_\bu(k+1,\bA) = \left\{(\zeta_1, \dots, \zeta_m)\in
\IR^m: \zeta_m \le \sum_{j=1}^{k+1} \lambda_j(A_m)\right\},$$
and every element in
$W_{k+1}(\bA) \cap \partial P_\bu(k+1,A)$ has the form
$\sum_{j=1}^{k+1} \la \bA x_j, x_j\ra$, where
$\{x_1, \dots, x_{k+1}\} \subseteq \cH$ is an
orthonormal family with
$(x_1, \dots, x_{r_1}) = (z_1, \dots, z_{r_1})$
and $A_m x_j = \xi_m x_j$ for $j > r_1$.

Consider
$$\bnu + \bxi \in \cl (\conv W_{k+1}(\bA))\cap \partial
P_\bu(k+1,\bA) = \conv W_{k+1}(\bA)\cap \partial P_\bu(k+1,\bA).$$
It is a convex combination of finitely many extreme points of
$\conv W_{k+1}(\bA)$ in
$\partial P_\bu(k+1,\bA)$. Since $W_{k+1}(\bA)$ is closed,
the extreme points of $\conv W_{k+1}(\bA)$ have the form
$\sum_{j=1}^{k+1} \la \bA x_j, x_j\ra$ with
$\sum_{j=1}^{k+1} \la A_m x_j, x_j\ra
= \sum_{j=1}^{k+1} \lambda_j(A_m)$.
We may assume that $(x_1, \dots, x_{r_1}) = (z_1, \dots, z_{r_1})$,
and
$x_j \in N(A_m - \xi_m I)$, the null space $A_m-\xi_m I$,
for $j = r_1 + 1, \ldots, k + 1$.
Hence, $N(A_m-\xi_m I)$
has dimension at least $k+1-r_1$.

Let $\cH_1 = \span(\{z_1, \dots, z_{r_1}\}\cup N(A_m-\xi_mI))$
and $\bB = (B_1, \dots, B_m)$, where $B_j$ is the compression of
$A_j$ onto $\cH_1$ for $j = 1, \dots, m$.

By the previous argument, for $\ell=k, k+1$, every element
in $W_\ell(\bA) \cap \partial P_\bu(\ell, \bA)$ has the form
$\sum_{j=1}^\ell \la \bA x_j, x_j\ra$ with
$(x_1, \dots, x_{r_1}) = (z_1, \dots, z_{r_1})$
and $x_j \in N(A_m - \xi_m I)$ for $j > r_1$.
As a result, for $\ell = k, k+1$,
$P_{\bu}(\ell,\bB) = P_\bu(\ell, \bA)$,
$$
W_{\ell}(\bB) \cap \partial P_\bu(\ell,\bB)
=
W_{\ell}(\bA)\cap \partial P_\bu(\ell,\bA),
$$
and
\begin{equation}\label{equal}
\conv W_{\ell}(\bB) \cap \partial P_\bu(\ell,\bB)
=
\conv W_{\ell}(\bA) \cap \partial P_\bu(\ell,\bA).
\end{equation}
If $\bnu \in \conv W_k(\bB)$, then
$\bnu \in \conv W_k(\bA)$. Since $\bnu$ is an extreme point of
$\cl (\conv W_k(\bA))$, it is also an extreme point of
$\conv W_k(\bA)$. Then we get $\bnu \in W_k(\bA)$, which
contradicts our assumption. So,
$$\bnu \notin \conv W_k(\bB)\cap
P_\bu(k,\bB) \subseteq \{(\zeta_1, \dots, \zeta_{m-1}, \nu_m):
\zeta_1, \dots, \zeta_{m-1}\in \IR\}.$$
Since $\bnu$ is an extreme point of
$\cl (\conv W_k(\bA))\cap \partial P_\bu(k,\bA)$,
it follows that $(\nu_1, \dots, \nu_{m-1})$
is an extreme point of the convex set
$$\{(\zeta_1, \dots, \zeta_{m-1}):
(\zeta_1, \dots, \zeta_{m-1}, \nu_m)
\in \cl (\conv W_k(\bA))\cap \partial P_\bu(k,\bA)\}.$$
So,
there is a
unit vector $\bv = (v_1, \dots, v_{m-1}) \in \IR^{m-1}$
such that
$$
\sum_{j=1}^{m-1} v_j \nu_j \ge \sum_{j=1}^{m-1} v_j \zeta_j$$
for any
$(\zeta_1, \dots, \zeta_{m-1}, \nu_m)
\in \cl (\conv W_k(\bA)) \cap \partial P_\bu(k,\bA)$,
which is the same as
$\cl (\conv W_k(\bB)) \cap \partial P_\bu(k,\bB)$
by (\ref{equal}).
We may apply the {\bf Reduction} argument to
$(B_1, \dots, B_{m-1})$ and assume that
$\bv = (0, \dots, 0, 1)$. As a result,
for any orthonormal family $\{x_1, \dots, x_k\}\subseteq \cH_1$
with
$(z_1, \dots, z_{r_1}) = (x_1, \dots, x_{r_1})$, we have
\begin{equation}
\label{nu(m-1)}
\nu_{m-1} \ge \sum_{j=1}^k \la A_{m-1} x_j, x_j\ra =
\sum_{j=1}^{k} \la B_{m-1} x_j, x_j\ra.
\end{equation}

Let
$\tilde \bmu = \sum_{j=1}^{r_1} \la \bA z_j, z_j\ra
= \sum_{j=1}^{r_1} \la \bB z_j, z_j\ra.$
Suppose $\hat \bA = (\hat A_1, \dots, \hat A_m)$
and $\hat \bB = (\hat B_1, \dots, \hat B_m)$
are such that for $j = 1, \dots, m$,
$\hat A_j$ is the compression of $A_j$ onto
$\{z_1, \dots, z_{r_1}\}^\perp$ in $\cH$, and
$\hat B_j$ is the compression of $B_j$ onto
$\{z_1, \dots, z_{r_1}\}^\perp$ in $\cH_1$.
 Then for $\ell = k, k+1$,
$$\cl (\conv W_\ell(\bA)) \cap P_\bu(\ell,\bA)
= \tilde \bmu + \cl (\conv W_{\ell-r_1}(\hat \bA))\cap P_\bu(\ell-r_1,\hat \bA)$$
and $W_\ell(\bB) \cap P_\bu(\ell, \bB) =
\tilde \bmu + W_{\ell-r_1}(\hat \bB)$.
We may do the following.

\medskip\noindent
{\bf Replacement}
\begin{itemize}
\item
Replace
$(\bA, \bB, k, k+1)$ by $(\hat \bA, \hat \bB, k-r_1, k+1-r_1)$,

\item Replace
$\cH$ by $\{z_1, \dots, z_{r_1}\}^{\perp}$ in $\cH$, and
replace
$\cH_1$ by $\{z_1, \dots, z_{r_1}\}^{\perp}$ in $\cH_1$.
\end{itemize}
After the replacement, we see that $B_m = \xi_m I_{\cH_1}$
and equation (\ref{nu(m-1)}) becomes
\begin{equation}\label{3.5}
\nu_{m-1} \ge \sum_{j=1}^k  \la A_{m-1}x_j, x_j\ra
= \sum_{j=1}^k  \la B_{m-1}x_j, x_j\rangle\end{equation}
for any orthonormal set $\{x_1, \dots, x_k\} \subseteq \cH_1$.
To carry out the next reduction, we need the following properties.
\begin{itemize}
\item[(I)] $\sum_{j=1}^r \la A_{m-1} z_j, z_j\ra =
\sum_{j=1}^r \la B_{m-1} z_j, z_j\ra
\ge \sum_{j=1}^r \lambda_j(B_{m-1})$, and

\item[(II)] $\lambda_r(B_{m-1}) \ge \xi_{m-1}$.
\end{itemize}

To prove (I),
let $\sum_{j=1}^{k+1} \la \bA y_j, y_j\ra$
be an extreme point of the non-empty closed convex set
$$\conv W_{k+1}(\bA)\cap P_\bu(k+1,\bA)$$
that
satisfies
\begin{equation}\label{3.6}
\sum_{j=1}^{k+1} \la A_{m-1} y_j, y_j \ra
= \max\{\zeta_{m-1}: (\zeta_1, \dots, \zeta_{m-1}, (k+1)\xi_m)
\in \conv W_{k+1}(\bA)\}, \end{equation}
where $\{y_1, \dots, y_{k+1}\}\subseteq \cH_1$ is an orthonormal
set.
Then since $\bnu+\bxi\in \conv W_{k+1}(\bA)\cap P_\bu(k+1,\bA)$,
$\sum_{j=1}^{k+1} \la A_{m-1} y_j, y_j \ra
\ge \nu_{m-1} + \xi_{m-1}.$
We have
\begin{equation}\label{3.7}
\la A_{m-1}y,y\ra \ge \xi_{m-1}
\quad \hbox{ for every unit vector }
y \in \span\{y_1, \dots, y_{k+1}\}.
\end{equation}
If (\ref{3.7}) does  not hold, we may apply a unitary transform to
$\{y_1, \dots, y_{k+1}\}$ and assume that $\la A_{m-1}y_{k+1},y_{k+1}\ra < \xi_{m-1}$.
Then
\begin{eqnarray*} \sum_{j=1}^{k} \la A_{m-1}y_j, y_j\ra & = &
\sum_{j=1}^{k+1} \la A_{m-1}y_j, y_j\ra
- \la A_{m-1}y_{k+1}, y_{k+1}\ra\\ & > &
(\nu_{m-1}+\xi_{m-1}) - \xi_{m-1} = \nu_{m-1},\end{eqnarray*}
contradicting (\ref{3.5}).
If (I) does not hold, there is an orthonormal set
$\{x_1, \dots, x_r\}$ in $\cH_1$ such that
$$ \sum_{j=1}^r \la A_{m-1} x_j, x_j\ra
= \sum_{j=1}^r \la B_{m-1} x_j, x_j\ra
> \sum_{j=1}^r \la B_{m-1} z_j, z_j\ra
= \sum_{j=1}^r \la A_{m-1} z_j, z_j\ra .$$
Observe that $\dim \big(\span\{y_1, \dots, y_{k+1}\} \cap \{x_1, \dots, x_{r}\}^\perp\big)\ge k-r + 1$.
By (\ref{3.7}), for any orthonormal set
$\{x_{r+1}, \dots, x_{k}\} \subseteq
\span\{y_1, \dots, y_{k+1}\} \cap \{x_1, \dots, x_{r}\}^\perp$,
we have $\la Ax_j, x_j \ra \ge \xi_{m-1}$ for $j = r+1, \dots, k$,
so that
$$
\sum_{j=1}^{k} \la A_{m-1} x_j, x_j \ra\\
> \sum_{j=1}^{r} \la A_{m-1} z_j, z_j \ra + \sum_{j=r+1}^{k}
\xi_{m-1} = \nu_{m-1},$$
contradicting (\ref{3.5}). Assertion (I) is proved.

To prove (II), apply a unitary transform to $\{z_1, \dots, z_r\}$
and assume that $B_{m-1}z_j = \lambda_j(B_{m-1}) z_j$ for
$j = 1, \dots, r$. If (II) does not hold, then
by (\ref{3.7}) we have
$$
\la B_{m-1} x_j, x_j \ra = \la A_{m-1} x_j, x_j\ra \ge \xi_{m-1},
\quad j = r, \dots, k,$$
for any orthonormal set
$\{x_r, \dots, x_k\} \subseteq \span\{y_1, \dots, y_{k+1}\}
\cap \{z_1, \dots, z_r\}^\perp$,
so that
\begin{eqnarray*} && \sum_{j=1}^{r-1} \la B_{m-1} z_j,z_j\ra +
\sum_{j=r}^k \la B_{m-1} x_j, x_j \ra\\ & > &
\left(\sum_{j=1}^{r} \la B_{m-1} z_j,z_j\ra - \xi_{m-1}\right)
+ (k-r+1) \xi_{m-1} = \nu_{m-1},\end{eqnarray*}
contradicting (\ref{3.5}).
Assertion (II) follows.

Let $r_2 = 0$ if $\lambda_1(B_{m-1}) =
\la B_{m-1}z_1, z_1\ra = \xi_{m-1}$, and
$r_2 = \max \{ j: \lambda_j(B) = \la B_{m-1} z_j, z_j \ra
> \xi_{m-1}\}$
otherwise. As in the previous step, the vectors $z_1, \ldots, z_{r_2}$
are uniquely determined up to a unitary transform. The null space $N(B_{m-1} - \xi_{m-1}I)$
is also large enough.
Indeed for the orthonormal set $\{y_1, \dots, y_{k+1}\}\subseteq \cH_1$ in the proof above,
we have
$$\sum_{j=1}^{k+1} \la B_{m-1} y_j, y_j \ra = \sum_{j=1}^{k+1}\lambda_j(B_{m-1})$$
so that the null space $N(B_{m-1} - \xi_{m-1}I)$ has dimension at least $k - r_2 + 1$.
Suppose $\cH_2$ is the subspace of $\cH_1$ spanned by
$\{z_1, \dots, z_{r_2}\} \cup N(B_{m-1} - \xi_{m-1}I)$.
Let $\bC = (C_1, \dots, C_m)$, where
$C_j$ is the compression of $A_j$ onto $\cH_2$ for $j = 1, \dots, m$.

If $\bnu \in \conv W_k(\bC)$, then $\bnu\in \conv W_k(\bA)$, and we get
a contradiction. So, $\bnu \notin \conv W_k(\bC)$ and is an extreme point of the closed convex set
$$\bS = \cl (\conv W_k(\bA))\cap
\left \{(\zeta_1, \dots, \zeta_{m}):
\zeta_{m-1} = \sum_{j=1}^k \lambda_j(B_{m-1}), \
\zeta_m =  k \xi_m \right\}.$$
Moreover, if $\{x_1, \dots, x_k\}$ is
an orthonormal set in $\cH_2$ such that
$\sum_{j=1}^k \la \bA x_j, x_j \ra \in \bS$,
then we have
$\sum_{j=1}^k \la B_{m-1} x_j, x_j\ra = \sum_{j=1}^k
\lambda_j(B_{m-1})$, and we may assume that
$(x_1, \dots, x_{r_2}) = (z_1, \dots, z_{r_2})$.
Thus, we may apply the
{\bf Reduction} and {\bf Replacement} steps
to $(C_1, \dots, C_m)$ and assume that
$(C_{m-1}, C_{m}) = (\xi_{m-1} I_{\cH_2}, \xi_{m}I_{\cH_2})$.
Then identify an orthonormal set
$\{\tilde y_1, \dots, \tilde y_{k+1}\} \subset \cH_2$
satisfying
$$\sum_{j=1}^{k+1} \la A_{m-2} \tilde y_j, \tilde y_j\rangle
= \max\{ \zeta_{m-2}:
(\zeta_1, \dots, \zeta_{m-2}, (k+1) \xi_{m-1}, (k+1)\xi_m)
\in W_{k+1}(\bA)\},$$
and repeat the arguments after (\ref{3.6})
with $(A_{m-1}, B_{m-1})$ replaced by $(A_{m-2}, C_{m-2})$.

Repeating the arguments, we will get
in $m-1$ steps $\bT
= (T_1, \dots, T_m) \in \cS(\cK)^m$
such that
$(T_2, \dots, T_m) = (\xi_2 I_\cK, \dots, \xi_m I_\cK)$,
where
$\cK = \cH_{m-1} \subseteq \cdots \subseteq \cH_1 \subseteq \cH$
is a nested sequence of subspaces  of $\cH$ such that
$$\nu_1 = \sum_{j=1}^r \la T_1 z_j, z_j \ra
+ (k-r) \xi_1,$$
and there is an orthonormal set
$\{\hat y_1, \dots, \hat y_{k+1}\}\subseteq\cH_{m-1}$ such that
\begin{eqnarray*}
&& \nu_1 + \xi_1
\sum_{j=1}^{k+1} \la T_1 \hat y_j, \hat y_j\ra\\ & = &
\max\{ \zeta_1: (\zeta_1, (k+1)\xi_2, \dots, (k+1)\xi_m) \in
W_{k+1}(\bT)\}\\
&=&\max\left \{ \sum_{j=1}^{k+1}\la T_1 x_j, x_j\ra:
\{x_1, \dots, x_{k+1}\} \subseteq \cH_{m-1}
\hbox{ is an orthonormal set} \right\}.
\end{eqnarray*}
Using the arguments after (\ref{3.6}), we see that
$\sum_{j=1}^{k+1}
\la T_1 \hat y_j, \hat y_j\ra = \sum_{j=1}^{k+1} \lambda_j(T_1)$.
So, we may apply
a unitary transform to $\hat y_1, \dots, \hat y_{k+1}$ and
assume that
$T_1 \hat y_j = \lambda_j(T_1) \hat y_j$ for $j = 1, \dots, k+1$,
and $\xi_1 = \lambda_{k+1}(T_1)$. As a result,
$\sum_{j=1}^k \la T_1 \hat y_j, \hat y_j\ra = \nu_1$
so that
$$
\bnu = \sum_{j=1}^k \la \bT \hat y_j, \hat y_j \ra =
\sum_{j=1}^k \la \bA \hat y_j, \hat y_j \ra \in W_k(\bA),$$
which is the desired contradiction. The proof is complete. \qed

\begin{remark} \label{4.8} \rm
It is tempting to use an inductive argument on $m$  in the proof of
Theorem \ref{k-closed} instead of doing the tedious
``Reduction and Replacement'' arguments. However,
as shown in the penultimate paragraph of the proof,
if $\bnu\in \conv W_k(\bC)$, one can quickly arrive at
the conclusion using an inductive argument. However,
if $\bnu \notin \conv W_k(\bC)$, more care is needed
in reducing the problem to the case with one fewer operator.
\end{remark}

In general, $W_k(\bA)$ is closed does not imply $W_{k+1}(\bA)$
is closed. For example, if
$A = I_k \oplus 0_k \oplus \diag(1, 1/2, 1/3,
\dots)$, then $W_k(A) = [0,k]$ and $W_{k+1}(A)  = (0,k+1]$.
It is interesting to prove or disprove that $W_k(\bA)$ is closed
whenever $W_{k+1}(\bA)$ is closed
for $\bA = (A_1, \dots, A_m)$ with $m> 3$.

Note that by Theorem \ref{k-closed}, if $W_{k+1}(\bA)$ is convex and closed,
then $\conv W_{k}(\bA)$ is closed.

\section{Polyhedral properties}
\setcounter{equation}{0}

In this section, we determine $\bA \in \BH^m$ such that $W_k(\bA)$ or
$\conv W_k(\bA)$ is a polyhedral set, i.e., the convex
hull of a finite set. It turns out that
these conditions are equivalent, and also equivalent
to that $A_1, \dots, A_m$ have
a common reducing subspace of finite dimension such that the compression
of $A_1, \dots, A_m$ onto the subspace are diagonal operators
$D_1, \dots, D_m$ and $W_k(A_1, \dots, A_m) = W_k(D_1, \dots, D_m)$.
This extends the result in \cite{LPW}.
It is remarkable that even if $\cH$ is inseparable, $W_k(\bA)$ is determined
by the compression of $A_1, \dots, A_m$ on a finite dimensional subspace if
$W_k(\bA)$ is polyhedral. We will obtain similar
descriptions for $\bA\in \BH^m$ such that $\cl(W_k(\bA))$ or
$\cl ( \conv W_k(\bA))$ is polyhedral.
We begin with the following.

\begin{proposition} \label{2.1} Suppose $\bA \in \BH^m$ and $k\in \IN$.
\begin{itemize}
\item[{\rm (a)}] The set $W_k(\bA)$ is a singleton if and only if
$A_j = \mu_j I$ for any $j = 1, \dots, m$.
\item[{\rm (b)}] The set $W_k(\bA)$ is a subset of a line segment in
$\IC^m$ if and only if there is $H \in \SH$ and complex numbers
$a_1, \dots, a_m, b_1, \dots, b_m$ such that
$A_j = a_j I + b_j H$ for $j = 1, \dots, m$.
\end{itemize}
\end{proposition}

\it Proof. \rm
(a)  The sufficiency is clear.
Conversely,  suppose $W_k(\bA)$ is a singleton.
If $A_1$ is not a scalar,
$W_k(A_1) = \{\mu: (\mu, \mu_2, \dots, \mu_m) \in W(\bA)\}$
is not a singleton, which is a contradiction.  Applying the same argument
to $A_j$ for any $j$, the result follows.

(b) The sufficiency is clear.  Conversely, suppose
$W_k(\bA)$ is a line segment. Then for any partial
isometry $X$ of finite rank,
$W_k(X^*A_1X, \dots, X^*A_mX)$ is a line segment
so that
$X^*A_jX = a_j(X) I + b_j(X) H(X)$
for some scalar $a_j(X), b_j(X)$ and Hermitian matrix
$H(X)$ depending on $X$
by \cite[Theorem 3.5]{LPW}.
Since this is true for any isometry, we conclude that
$A_1, \dots, A_m$ have the said structure.
\qed

To facilitate our discussion, we present some definitions and notations.
To describe the connections among
the operators $A_1, \dots, A_m$, we often say that there is a unitary
operator $U = [U_1 \ U_2] \in \BH$ such that for each $j = 1, \dots, m$,
$U^*A_jU = \begin{pmatrix} A_{11}^{(j)} &  A_{12}^{(j)} \cr
 A_{21}^{(j)} &  A_{22}^{(j)}\cr\end{pmatrix}$,
where $A_{rs} =  U_r^*A_j U_s$ for $r,s \in \{1,2\}$ with some
special properties. For example, denote by $\cH_1$ the range space of
the partial isometry $U_1$.

\begin{enumerate}
\item  The operator
 $A_{21}^{(j)}$ is zero for each $j$
 means that $\cH_1$ is a common invariant subspace
 of $A_1, \dots, A_m$.
\item The operators $A_{12}^{(j)}$ and $A_{21}^{(j)}$ are zero for all $j$, i.e.,
$A_j = A_{11}^{(j)} \oplus A_{22}^{(j)}$, means that
$\cH_1$ is a common reducing subspace of $A_1, \dots, A_m$.
\item The $m$-tuple of
operators $(A_{11}^{(1)}, \dots, A_{11}^{(m)})$ is a
{\it (common) compression} of $\bA = (A_1, \dots, A_m)$ onto the
subspace $\cH_1$.
\end{enumerate}

Let $\cS$ be a convex subset of $ \IR^m$ (or $\IC^m$
identified as $\IR^{2m}$) and $\bx\in \cS$.
With respect to the standard inner product $\la\ ,\ \ra$, a nonzero vector
$\bv\in \IR^m$ defines a hyperplane
$$P(\bv,\bx)=\{\bu\in \IR^m:\langle \bv,\bu -\bx\rangle =0\}$$
and two half spaces $H_+(\bv,\bx)=\{\bu\in \IR^m:\langle \bv,\bu -\bx\rangle \ge0\}$ and $H_-(\bv,\bx)=\{\bu\in \IR^m:\langle \bv,\bu -\bx\rangle \le0\}$. $\bv$ is called a normal vector of $P(\bv,\bx)$.
The hyperplane $P(\bv,\bx)$ is called a {\it supporting plane} of $\cS$ at $\bx$ if
$\cS$ is a subset of either
$H_+(\bv,\bx)$ or $H_-(\bv,\bx)$.
A boundary point
$\bx$ of $\cS$ is
a {\it conical point} of $\cS$ if there are $m$ supporting planes of
$\cS$ at $\bx$ with linearly
independent normal vectors.
Clearly every extreme point of a polyhedral set is a conical point.

\begin{theorem} \label{2.3}
Let $\bA = (A_1, \dots, A_m)\in \BH^m$ and $k \in \IN$.
The following conditions are equivalent.
\begin{itemize}
\item[{\rm (a)}] $W_\ell(\bA)$ is polyhedral for all
$\ell = 1, \dots, k$.
\item[{\rm (b)}] $\conv W_k(\bA)$ is polyhedral.
\item[{\rm (c)}] There is $r \in \IN$ with $r\ge 2k$ and a unitary $U \in \BH$
such that for each $j = 1, \dots, m$,  $U^*A_jU = D_j \oplus B_j$, where
$D_j \in M_r$ is a diagonal matrix, satisfying
$W_k(\bA) = W_k(D_1, \dots, D_m)$.
\end{itemize}
\end{theorem}

\it Proof. \rm  We may assume that $\bA \in \SH^m$.
The implication (a) $\Rightarrow$ (b) is clear.

Suppose (b) holds so that $\conv W_k(\bA)$ is polyhedral
with $q$ vertices $\bnu_1, \dots, \bnu_q \in \IR^m$.
For each $p = 1, \dots, q$,
there is an orthonormal set
$S_p = \{x_1^{(p)}, \dots, x_k^{(p)}\}$ such that
$\bnu_p =
\sum_{j=1}^k \la \bA x_j^{(p)}, x_j^{p)}\ra$.
Let $V$ be any finite dimensional subspace of $\cH$ containing $\span S_p$ (for all $p = 1, \ldots, q$), and
$\tilde A_j$ be the compression of $A_j$ onto $V$ for
$j = 1, \dots, m$. Then for
$\tilde \bA = (\tilde A_1, \dots, \tilde A_m)$,
$\bnu_p \in
\conv W_k(\tilde \bA) \subseteq \conv W_k(\bA)$ so that $\bnu_p$ is a conical point of $\conv W_k(\tilde \bA)$.
By \cite[Theorem 4.2]{LPW}, $\span S_p$ is a common reducing subspace for $\tilde A_1, \dots, \tilde A_m$.
As this is true for any such subspace $V$, $\span S_p$ is a common reducing subspace
for $A_1, \ldots, A_m$. Let $W = \span(S_1 \cup \cdots \cup S_q)$. It is clearly an invariant subspace for all $A_1, \ldots, A_m$.
But as they are self-adjoint, $W$ is indeed a common reducing subspace. Now let $\tilde A_j$ be the compression of $A_j$ onto $W$ for $j = 1, \dots, m$. We have
$$\conv\{\bnu_1, \dots, \bnu_q\}
= \conv W_k(\tilde A_1, \dots, \tilde A_m)
= \conv W_k(\bA).$$
By \cite[Theorem 5.1]{LPW},
there is a unitary $\tilde U$ on $W$ such that
$\tilde U^*\tilde A_j \tilde U = D_j \oplus B_j$,
where $D_j \in M_{r_1}$ is a diagonal matrix
for $j = 1,\dots, m$, and
$$W_k(D_1,\dots, D_m) = W_k(\tilde A_1, \dots, \tilde A_m) =
\conv W_k(\tilde A_1, \dots, \tilde A_m).$$
The operator $U = \tilde U\oplus I_{W^\perp}$ is unitary and satisfies our purpose.

Next, we show that $D_j\in M_r$ can be chosen with $r \ge 2k$. To this end, let
$\tilde V$ be a subspace of dimension $2k+r_1$ containing $W$.
Suppose $C_j$ is the compression of
$A_j$ to $\tilde V$ for $j = 1,\dots, m$.
Then for $\bC = (C_1, \dots, C_m)$ and $\bD = (D_1,\dots, D_m)$,
$$W_k(\bD) \subseteq W_k(\bC)
\subseteq W_k(\bA) = W_k(\bD).$$
By \cite[Theorem 5.1]{LPW} again, we can choose   $D_j\in M_r$ with $r \ge 2k$.  Hence (b) $\Rightarrow$ (c).

Now, suppose (c) holds.
For every unit vector $\bu =(u_1, \dots, u_m) \in \IR^m$,
\begin{eqnarray*}
&& \max\{\bu\cdot \bxi: \bxi\in W_k(\bD)\}\\
& = &
\mbox{$\max\left\{\sum_{j = 1}^m u_j \sum_{\ell = 1}^k \la D_j x_\ell, x_\ell\ra: \{x_1, \ldots, x_k\}\ \mbox{is an
o.n. set in $\cH$}\right\}$}\\
& = &
\mbox{$\max\left\{\sum_{\ell = 1}^k \left\langle\left(\sum_{j = 1}^m u_jD_j\right) x_\ell, x_\ell\right\rangle: \{x_1, \ldots, x_k\}\ \mbox{is an
o.n. set in $\cH$}\right\}$}
\end{eqnarray*}
is the sum of the
$k$ largest eigenvalues of $\bD_\bu = \sum_{j=1}^m u_j D_j$.
Similarly,
if $\bA_\bu = \sum_{j=1}^m u_j A_j$, then
$\max\{\bu\cdot \bxi: \bxi\in W_k(\bA)\}
= \sum_{j = 1}^k \lambda_j(\bA_\bu)$.
So, $W_k(\bA) = W_k(\bD)$ implies that $\sum_{j = 1}^k \lambda_j(\bA_\bu) = \sum_{j = 1}^k \lambda_j(\bD_\bu)$.
Since $\bA_\bu$ is unitarily similar to a direct sum of the form $\bD_\bu \oplus \bB_\bu$,
$\lambda_j(\bA_\bu)\ge \lambda_j(\bD_\bu)$
 for $j = 1, \ldots, k$.
We see that $\lambda_j(\bA_\bu) = \lambda_j(\bD_\bu)$
is actually the $j$th largest eigenvalue of $\bA_\bu$ for each $j$.
It follows that for every $\ell = 1, \dots, k$,
$$\max\left\{\bu\cdot \bxi: \bxi\in W_\ell(\bA)\right\}
= \sum_{j = 1}^\ell \lambda_j(\bA_\bu) = \sum_{j = 1}^\ell \lambda_j(\bD_\bu) =
\max\left\{\bu\cdot \bxi: \bxi\in W_\ell(\bD)\right\}.$$
As this is true for all unit vectors
$\bu$, $\conv W_\ell(\bA)$ and $W_\ell(\bD)$
have the same support plane in all directions. Also, as $\conv W_\ell(\bA)$ is closed
by Theorem \ref{k-closed},
the two sets are equal.
Thus,
$$W_\ell(\bD) \subseteq W_\ell(\bA)
\subseteq \conv W_\ell(\bA)  = W_\ell(\bD),$$
and so, $W_\ell(\bA) = W_\ell(\bD)$ is polyhedral. We have (c) $\Rightarrow$ (a), and the proof is complete.
\qed

\medskip
We have an approximation result without assuming that
$W_k(\bA)$ or $\conv W_k(\bA)$ is closed.

\begin{theorem} \label{2.4} Let $\bA = (A_1, \dots, A_m)\in \BH^m$
and $k \in \IN$.
The following conditions are equivalent.
\begin{itemize}
\item[{\rm (a)}] $\cl W_\ell(\bA)$ is polyhedral
for all $\ell = 1, \dots, k$.
\item[{\rm (b)}] $\cl (\conv W_k(\bA))$
 is polyhedral.
\item[{\rm (c)}] There exists $r \ge 2k$ and for each $n\in \IN$ there is
an isometry
$X_n: \IC^r \rightarrow \cH$ such that
$\bD_n = (X_n^*A_1X_n, \dots, X_n^*A_m X_n)$
is an $m$-tuple of $r\times r$ diagonal matrices
satisfying
$$\bD_n \rightarrow \bD \quad \hbox{ and } \quad
W_k(\bD_n) \rightarrow W_k(\bD) = \cl(W_k(\bA)),$$
relative to the Hausdorff distance.
 \end{itemize}
\end{theorem}

\it Proof. \rm We may assume that $\bA \in \SH^m$.
The implication (a) $\Rightarrow$ (b) is clear.

To prove (b) $\Rightarrow$ (c), suppose (b) holds so that  $\cl(\conv W_k(\bA))$ is polyhedral
with $q$ vertices $\bnu_1, \dots, \bnu_q$.
By Theorem \ref{closure}, $\bnu_j = \bmu_j + (k-\ell_j)\bxi_j$
with $\ell_j\in \{0, \dots, k\}$, $\bmu_j \in W_{\ell_j}(\bA)$
and $\bxi_j \in W_{\rm ess}(\bA)$.
For $\bv = (v_1, \dots, v_m) \in \IR^m$, and $p \in \IN$,
let $\bv I_p = (v_1 I_p,\dots, v_m I_p)$.
Set $\tilde \cH = \IC^p \oplus \cH$ for
$p = mk - \ell_1 - \cdots - \ell_p$,
$$\bQ = (Q_1, \dots, Q_m) = \bxi_1 I_{k-\ell_1} \oplus
\cdots \oplus \bxi_q I_{k-\ell_q} \quad \hbox{ and }
\quad\tilde \bA = \bQ \oplus \bA \in \mathcal S(\tilde \cH)^m.$$
Then clearly, $\cl (\conv W_k(\bA)) \subseteq \cl (\conv W_k(\tilde \bA))$.
On the other hand, it is not difficult to adapt the finite dimensional result \cite[Theorem 4.4]{LPW}
to get
\begin{eqnarray*}
W_k(\tilde \bA) & \subseteq & \conv\bigcup_{\ell = 0}^k \big[W_{\ell}(\bQ) + W_{k - \ell}(\bA)\big]\\
& \subseteq & \conv\bigcup_{\ell = 0}^k \big[\ell W_{\rm ess}(\bA) + W_{k - \ell}(\bA)\big]
\ \subseteq\ \cl (\conv W_k(\bA)).
\end{eqnarray*}
So, $\cl (\conv W_k(\bA)) = \cl (\conv W_k(\tilde \bA))$.
By construction,
$W_k(\tilde \bA)$ contains all the vertices of $\cl (\conv W_k(\tilde \bA))$.
It follows that
$\conv W_k(\tilde \bA) = \cl (\conv W_k(\bA))$ is polyhedral.
By Theorem  \ref{2.3},
$\tilde \bA$ is unitarily similar to
$\bD \oplus \bB$, where $\bD = (D_1, \dots, D_m)$
is an $m$-tuple of diagonal matrices
in $M_r$ with $r \ge 2k$, and
$W_k(\tilde \bA) = W_k(\bD)$.
Suppose $D_j = \diag(\xi_{j1}, \dots, \xi_{jr})$.
There is an orthonormal set
$\{u_1, \dots, u_r\}$ such that
$\tilde A_j u_\ell = \xi_{j\ell} u_\ell$,
for each $j = 1, \dots, m$ and $\ell = 1, \dots, r$.
Let $u_\ell = \begin{pmatrix} v_\ell\cr w_\ell\cr\end{pmatrix}$
according to the decomposition $\tilde \cH = \IC^p \oplus \cH$.
Then $D_j v_\ell  = \xi_{j\ell} v_\ell$
and $A_j w_\ell = \xi_{j\ell} w_\ell$ for all $j$ and $\ell$.
Now, decompose $\tilde \cH$ as $\span\{v_1, \dots, v_r\}
\oplus \span\{w_1, \dots, w_r\} \oplus \tilde H_0$,
where $\tilde H_0 = \{v_1, \dots, v_r, w_1, \dots, w_r\}^\perp$
in $\tilde \cH$. Using this space decomposition, we may
assume that
$$\tilde \bA = \tilde \bD_Q \oplus \tilde \bD_A \oplus \tilde \bB$$
such that $W_k(\tilde \bA) = W_k(\tilde \bD_Q \oplus \tilde \bD_A)$,
where $\tilde \bD_Q$ is the compression of $\tilde \bA$ onto
$\span\{v_1, \dots, v_r\}$,
$\tilde \bD_Q$ is the compression of $\tilde \bA$ onto
$\span\{w_1, \dots, w_r\}$, $\tilde \bD_Q \oplus \tilde \bD_A$
is an $m$-tuple of diagonal matrices in $M_{\tilde r}$
with $\tilde r \ge r$,
and $\tilde \bB$ is the compression of $\tilde \bA$ onto $\cH_0$.
Actually, there is no harm to assume that $\tilde \bD_Q = \bQ$,
$\bA = \tilde \bD_A \oplus \tilde \bB$
with $\tilde \bB = (\tilde B_1, \dots, \tilde B_m)$,
and we may reset $r = \tilde r$.

Let $n\in\mathbb N$, and  let
$\{y_1, \ldots, y_{r_1}\}$ and $\{x_1, \ldots, x_{r_2}\}$ be
complete sets of orthonormal eigenvectors of $\bQ$ and $\tilde \bD_A$ respectively.
Since $\la \bQ y_1, y_1\ra \in W_{\rm ess}(\bA)$,
$\la \bQ y_1, y_1\ra \in W_{\rm ess}(\tilde \bB)$ as well.
We can find a unit vector $x_{r_2 + 1}\in\{x_1, \ldots, x_{r_2}\}^\perp$
such that
$$\|\la \bA x_{r_2 + 1}, x_{r_2 + 1}\ra - \la \bQ y_1, y_1\ra \|_2
= \|\la \tilde \bB x_{r_2 + 1}, x_{r_2 + 1}\ra - \la \bQ y_1, y_1\ra\|_2 < 1/n.$$
Here $\|\cdot\|_2$ denotes the usual Euclidean norm in $\mathbb R^m$.
Let $\tilde \bB_1$ be the compression of $\bA$ onto
$\cH_1 = \{x_1, \ldots, x_{r_2}, x_{r_2 + 1}, A_1x_{r_2 + 1}, \ldots, A_mx_{r_2 + 1}\}^\perp$.
Then $\la\bQ y_2, y_2\ra\in W_{\rm ess}(\tilde \bB_1)$. So,
there is a unit vector $x_{r_2 + 2}\in \cH_1$
such that $\|\la \bA x_{r_2 + 2}, x_{r_2 + 2}\ra - \la \bQ y_2, y_2\ra \|_2 < 1/n$.
Note that the compression of each $A_j$ onto $\span\{x_1, \dots, x_{r_2 + 1}, x_{r_2 + 2}\}$ is diagonal.
Inductively, we are able to find an orthonormal set $\{x_1, \dots, x_{r_2 + r_1}\}$
such that $\|\la \bA x_{r_2 + \ell}, x_{r_2 + \ell}\ra - \la \bQ y_\ell, y_\ell \ra \|_2 < 1/n$
for all $\ell = 1, \ldots, r_1$, and the compression of each $A_j$ onto
$\span\{x_1, \dots, x_{r_2 + r_1}\}$ is diagonal. Denote the compression of $\bA$ onto
$\span\{x_1, \dots, x_{r_2 + r_1}\}$ as $\bD_n$.

For each vertex $\bmu_p$ of $\cl(\conv W_k(\bA))
= W_k(\bQ\oplus \tilde \bD_A)$,
$$\bmu_p = \sum_{j = 1}^\ell \la \bQ z_j, z_j\ra + \sum_{j = \ell + 1}^{k - \ell} \tilde \bD_A z_j,$$
where $z_1, \ldots, z_\ell\in \{y_1, \ldots, y_{r_1}\}$ and $z_{\ell + 1}, \ldots, z_k\in \{x_1, \ldots, x_{r_2}\}$.
By suitably replacing those $z_1, \ldots, z_\ell\in \{y_1, \ldots, y_{r_1}\}$,
it is easy to find a $\bnu_p\in W_k(\bD_n)$ such that
$\|\bmu_p - \bnu_p\|_2 < k/n$.
Now each $\bmu\in \conv W_k(\bA)$ can be written as a convex combination
$\sum_{p = 1}^q \lambda_p\bmu_p$ of the vertices. Then $\bnu = \sum_{p = 1}^q \lambda_p\bnu_p$
is a point in $W_k(\bD_n)$ such that $\|\bmu - \bnu\|_2 < k/n$. As $\bD_n$ is a compression
of $\bA$, $W_k(\bD_n)\subseteq W_k(\bA)\subseteq \cl(\conv W_k(\bA))$.
The Hausdorff distance between the compact sets
$W_k(\bD_n)$ and $\cl(\conv W_k(\bA))$ is no more than $k/n$.
So, $W_k(\bD_n) \rightarrow \cl(\conv W_k(\bA))$.  Moreover,
$$\cl(\conv W_k(\bA))\subseteq \cl(\cup_{n = 1}^\infty W_k(\bD_n))
\subseteq \cl W_k(\bA)\subseteq \cl(\conv W_k(\bA)),$$
and we get $\cl W_k(\bA) = \cl(\conv W_k(\bA))$.
Finally, identifying $\bD_n$ with its matrix representation,
then by considering a subsequence, we may
assume that $\bD_n = (D_1^{(n)}, \dots, D_m^{(n)})$ converge
to $\bD = (D_1, \dots, D_m)$.
Since we also have $W_k(\bD_n)\to W_k(\bD)$, $\cl W_k(\bA)= W_k(\bD)$.
Condition (c) follows.

Now, if (c) holds, then $\cl W_k(\bA)= W_k(\bD)$
is polyhedral. Condition (a) holds with $\ell = k$.
\medskip
For $\ell < k$, write
for every unit vector $\bu = (u_1, \dots, u_m) \in
\IR^m$, $A_\bu  = \sum_{j=1}^m u_j A_j$
and $D_\bu = \sum_{j=1}^m u_j D_j$. Since
$\cl W_k(\bA)= W_k(\bD)$,
we get as in the proof of Theorem  \ref{2.3},
$\sum_{j=1}^k \lambda_j(A_\bu) = \sum_{j=1}^k \lambda_j(D_\bu)$.
Also, as $D_\bu$ is the limit of a sequence of principal submatrices
of $A_\bu$, $\lambda_j(D_\bu)\le \lambda_j(A_\bu)$ for $j =1 , \dots, k$,
by Cauchy's interlacing theorem.
Consequently, if $\ell \in\{1, \dots, k\}$, then
$\sum_{j=1}^\ell \lambda_j(A_\bu) = \sum_{j=1}^\ell \lambda_j(D_\bu)$.
Since this is true for every unit vector $\bu$,
$\cl(\conv W_\ell(\bA)) = W_\ell(\bD)$.
On the other hand, $W_\ell(\bD) \subseteq \cl W_\ell(\bA)$ as $\bD$ is a limit of an $m$-tuple of
principal submatrices of $\bA$.
As above,
$$W_\ell(\bD) \subseteq \cl W_\ell(\bA)
\subseteq \cl(\conv W_\ell(\bA))  = W_\ell(\bD).$$
Hence, $\cl W_\ell(\bA) = W_\ell(\bD)$ is polyhedral,
and condition (a) holds.
\qed

One readily see that every vertex $(\mu_1, \ldots, \mu_m)$
of $W_k(\bA)$, $\mu_j$ is the sum of $k$ approximate eigenvalue of
$A_j$. Using the ideas in the proof of Theorems \ref{2.3} and \ref{2.4},
we have the following result. Part (a) is an extension
of the result in \cite{MF} for a single matrix, and part (b) is an extension of
the result in \cite{Cho2}.

\begin{theorem} \label{5.4}
Let $\bA \in \BH^m$ and $k \in \IN$.
\begin{itemize}
\item[{\rm (a)}] If $\bmu = (\mu_1, \dots, \mu_m)
= (\tr(X^*A_1X), \dots, \tr(X^*A_mX))$ is a conical point of
$W_k(\bA)$, where
$X: \IC^k \rightarrow \cH$ is an isometry,
then the range space of $X$ is a reducing subspace of $A_j$ for $j = 1, \dots, m$,
i.e., $U^*A_jU = C_j \oplus B_j$ for $j =1, \dots, m$,  if $U = [X | X^\perp]$.
\item[{\rm (b)}] If $\bmu = (\mu_1, \dots, \mu_m)$
is a conical point of $\cl(W_k(\bA))$ and not in $W_k(\bA)$, then there
is a sequence of isometries $X_\ell: \IC^k\rightarrow \cH$
for $\ell \in \IN$
such that
$$\lim_{\ell\rightarrow \infty}(\tr(X_\ell^* A_1 X_\ell), \dots,
\tr(X_\ell^* A_m X_\ell)) =\bmu.$$
\end{itemize}
\end{theorem}

\it Proof. \rm We may assume that $A_1,\dots, A_m$ are self-adjoint.
(a) Suppose $\bmu = \sum_{j=1}^k \la \bA x_j, x_j\ra  \in W_k(\bA)$
is a conical point, where $\{x_1, \dots, x_k\}$ is an orthonormal set.
Consider the operator matrices $A_j =
\begin{pmatrix} C_j & Z_j \cr Z_j^* & B_j\cr\end{pmatrix}$
for $j = 1, \dots, m$,
where $C_j \in M_k$ is the compression of $A_j$ onto the space
spanned by $\{x_1, \dots, x_k\}$. We will show
that $Z_j = 0$ for all $j$ and the result will follow.
Assume without loss of generality that $Z_1 \ne 0$, and
by abuse of notation that the $(1,1)$ entry of $Z_1$ is nonzero.
Consider the leading principal submatrix $T \in M_{k+1}$
of $A_1 + iA_2$. Then
$$\mu_1 + i \mu_2 \in W(T) \subseteq
\{\xi_1 + i\xi_2: (\xi_1, \dots, \xi_m) \in W(\bA)\}$$
so that $\mu_1 + i\mu_2$ is a conical point of $W(T)$.
By the result in \cite{MF}, $T = T_1 \oplus T_2$ with $T_1 \in M_k$,
which is a contradiction.

(b) Suppose $\bmu$ is a conical point of $\cl W_k(\bA)$.
Then $\bmu$ is an extreme point of $\cl (\conv W_k(\bA))$.
Thus, $\bmu = \bnu + (k-\ell) \bxi$ with
$\bnu \in W_\ell(\bA)$ and $\bxi \in W_{\ess}(\bA)$ for some
$\ell \in \{0, \dots, k\}$. Let $\tilde \bA = \bA\oplus \bxi I_{k-\ell}$.
Then $\bmu$ is a conical point of $W_k(\tilde \bA)$.
Using (a) and an argument similar to the proof of Theorem \ref{2.4},
we get the conclusion. \qed

\section{Compact operators}
\setcounter{equation}{0}

In \cite[Theorem 3.3]{LPW}, it was shown that $\{A_1, \dots, A_m\} \subseteq M_n$
is a commuting family of normal matrices if and only if $W_k(A_1, \dots, A_m)$
is poyhedral for all $k = 1, \dots, n-1$. In this section, we study conditions
for compact operators $A_1, \ldots, A_m$ satisfying one of the
following conditions.

\begin{itemize}
\item[{\rm (a)}]  $W_k(\bA)$ is polyhedral for all $k \in \IN$.

\item[{\rm (b)}] $A_1, \ldots, A_m$ are commuting normal operators.
\end{itemize}
One may see some results for the joint numerical range of compact operators
in \cite{Barra,Juneja} and their references.
There have been considerable interests
in the commutativity of compact self-adjoint (normal) operators in
terms of different kinds of joint spectra. See for example \cite{Cho},
\cite{GR} and \cite{MQW}. Our result gives an alternative
from the joint numerical range viewpoint, and use the joint eigenvalues
of operators defined as follows.

Let $\bA \in \BH^m$. A point $\bmu = (\mu_1, \ldots, \mu_m)\in \mathbb C^m$
is a {\it joint eigenvalue}
of $\bA$ if there is a nonzero $x\in \cH$ such that
$A_j x = \mu_j x$ for all $j = 1, \ldots, m$. The nonzero $x$ is a common eigenvector corresponding to $\bmu$.
Define $\Sigma_k(\bA)$ to be the set of all $\bmu_1 + \cdots + \bmu_k$, where $\bmu_1, \ldots, \bmu_k$ are joint eigenvalues of $\bA$ corresponding to $k$ linearly independent common eigenvectors $x_1, \dots, x_k$.
This set could be empty; indeed, $\bA$ may not have any joint eigenvalue at all.
When $k = 1$, $\Sigma_1(\bA)$ reduces to the joint
point spectrum of $\bA$.
Suppose $\bmu_1, \ldots, \bmu_k$ are joint eigenvalues of $\bA$ corresponding
to linearly independent common eigenvectors $x_1, \dots, x_k$.
One may apply Gram-Schmidt process to $\{x_1, \dots, x_k\}$ to obtain an orthonormal
set $\{u_1, \dots, u_k\}$. Then the matrix
$B_j = ( \la A_j u_r, u_s\rangle )\in M_k$ is in
upper triangular form, and
$\bmu_\ell = ((B_1)_{\ell\ell}, \dots, (B_m)_{\ell\ell})$
for $\ell = 1, \dots, k$.
Thus,  $\la \bA u_\ell, u_\ell\ra = \bmu_\ell$ for $\ell = 1, \dots, k$
and $\sum_{j=1}^k \bmu_j
= \sum_{j=1}^m \la \bA u_j, u_j\ra \in  W_k(\bA)$. As a result,
 $\Sigma_k(\bA)\subseteq
W_k(\bA)$.
Evidently, conditions (a) -- (c) in Theorem \ref{2.3} are
equivalent to the condition:

\medskip
(d) \it  $W_k(\bA) = \conv\Sigma_k(A)$ and\, $\conv \Sigma_k(A)$ is polyhedral.

\rm
\medskip\noindent
Suppose $A_1, \dots, A_m$ are commuting compact
normal operators in $\BH$ and $\bA = (A_1, \dots, A_m)$.
As far as its $k$-numerical ranges are concerned, we may assume that $\cH$ is separable.
By \cite[Proposition 1]{GR}, there is a unitary $U$ such that
$$D_j=U^*A_jU=\diag(d_j(1), d_j(2),\dots )$$
is a diagonal operator for $j = 1,\dots, m$.
Let $\bv(\ell)=(d_1(\ell), d_2(\ell),\dots , d_m(\ell))$ for  $\ell\ge 1$. Then
for $k\in \IN$,
$$\Sigma_k(\bA)=\left\{\sum_{\ell=1}^k\bv(j_\ell):1\le j_1<j_2<\cdots<j_k\right\}.$$
For $k\ge 1$, let $\cP_k$ be the set of rank $k$
orthogonal projections in $\BH$, and
$\cD_k=\{\diag P:P \in \cP_k\}$. By \cite[Theorem 13]{Kadison},
$\cD_k$ consists of all
$(u_1,u_2,\dots)$ such that
$
0\le u_j\le 1$ for $j=1,\dots$ and
$\sum_{j=1}^{\infty}u_j=k$. Then we have
  $$W_k(\bA)
  =
  \left\{\(\sum_{j=1}^{\infty}u_jd_1(j), \sum_{j=1}^{\infty}u_jd_2(j),
  \dots , \sum_{j=1}^{\infty}u_jd_m(j)\):(u_1,u_2,\dots)\in \cD_k
  \right\}.$$
Suppose $\dim \cH=N$ is finite and $1\le k\le N$.
Let $\bv=(\underbrace{1,\dots,1}_{k-{\rm terms}},0,\dots,0)\in \IR^N$.
Then,
$$\cD_k
=\conv\{P\bv:P\mbox{ is an }N\times N\mbox{ permutation matrix}\},$$
for example see \cite{MO}.
Therefore, we have
\begin{equation}\label{sigma}
W_k(\bA)= \conv \Sigma_k(\bA) .
\end{equation}
In fact, it is known that $W_k(\bA) = \conv \Sigma_k(\bA)$ for all
$k = 1, \dots, N-1$ if and only if $\{A_1, \dots, A_m\}$
is a commuting family; see \cite{LPW} and also \cite{LST,MMF}.
The result can be extended to compact operators
as in the following.

\begin{theorem} \label{3.1}
Suppose $\bA=(A_1, \dots, A_m)\in \BH^m$ is
an $m$-tuple of compact operators.
Then  $\{A_1, \dots, A_m\}$ is a commuting family of normal operators
if and only if
$$W_k(\bA)= \conv \Sigma_k(\bA)\qquad \hbox{ for all } \ k \ge 1.$$
\end{theorem}

\it Proof. \rm Suppose $\{A_1, \dots, A_m\}$ is a commuting family of compact normal operators.
Without loss of generality, we may assume that $\cH$ is separable,
$A_j=\diag(d_j(1), d_j(2),\dots )$  for $j=1,\dots,m$ and
$\bv(\ell)=(d_1(\ell), \dots , d_m(\ell))\in \IC^m$
for all $\ell\ge 1$.  Let $\bmu = \sum_{j = 1}^\infty u_j\bv(j)$ be a point in $W_k(\bA)$, where $(u_1,u_2,\dots)\in \cD_k$.

If $u_j\neq 0$ for only finitely many $j$'s, then by (\ref{sigma}), $\bmu\in \conv\Sigma_k(\bA)$. Otherwise,
choose $N > k$ such that $a=\sum_{j =N+ 1}^{\infty} u_j< 1$
and consider the sum
$\bu = \sum_{j = N + 1}^\infty (u_j/a)\bv(j)$.
It is an infinite convex combination of the points $\bv(j)$, $j = N + 1, \ldots$ in $\mathbb R^m$.
By the result of \cite{CW}, there exists $\ell >0$ such that $\bu = \sum_{i= 1}^\ell \alpha_i\bv(N + i)$ for some
$\alpha_i\ge 0$ satisfying $\sum_{i= 1}^\ell \alpha_i = 1$. We have
$$\bmu = \sum_{j = 1}^N u_j\bv(j)+a\bu
= \sum_{j = 1}^N u_j\bv(j)+ \sum_{i = 1}^\ell a\alpha_i\bv(N + i)$$
with
$$0\le u_j,\, a\alpha_i\le 1 \ \mbox{ for all }1\le j\le N,\ 1\le i\le \ell \quad \mbox{ and }
\quad \sum_{j = 1}^N u_j\ + \sum_{i = 1}^\ell  a\alpha_i=k.$$
By (\ref{sigma}) again, $\bmu \in \conv\Sigma_k(\bA)$.

Conversely, suppose $W_k(\bA) = \conv \Sigma_k(\bA)$ for all
$k \ge 1$.
Consider the first operator $A_1$. We have $W_k(A_1) = \conv \Sigma_k(A_1)$ for all $k$.
Write $A_1 = H_1 + iG_1$ for self-adjoint $H_1, G_1\in \BH$.
The compact self-adjoint $H_1$ can be written as
$$H_1 = \sum_{j \ge 1} a_j I_{\cH_j} + \sum_{\ell \ge 1} b_\ell I_{\cK_\ell},$$
where $a_j$'s are the positive, and $b_\ell$'s the negative, eigenvalues of $H_1$.
The eigenspaces $\cH_j$'s and $\cK_\ell$'s are mutually orthogonal with $\dim \cH_j = n_j$ and
$\dim \cK_\ell = m_\ell$ all finite.
Either sum can be empty, finite, or infinite. An empty sum is taken to be 0.
We shall also assume that the eigenvalues are ordered as
$a_1 > a_2 > \cdots$ and $b_1 < b_2 < \cdots$.
We shall show that
$\cH_j$'s and $\cK_\ell$'s are all reducing subspaces of $G_1$. It follows that
the null space of $H_1$ is also a reducing subspace of $G_1$, and consequently, $H_1$ and $G_1$ commute,
and $A_1$ is normal.

Starting with $a_1$, let $k = \dim \cH_1 = n_1$.
Clearly, $\max \re W_k(A_1) = n_1a_1$.
Take an extreme point $\mu$ of $\cl W_k(A_1)$ lying on the vertical line $\re z = n_1a_1$.
Then $\mu\in W_k(A_1)$. For otherwise $\mu = \mu_1 + (k - \ell) 0 = \mu_1$ for $\mu_1\in W_{\ell}(A_1)$ with
$\ell < k$ (Theorem \ref{closure}),
and so $\re \mu = \re \mu_1 < n_1a_1$, a contradiction.
Hence $\mu$ is an extreme point of $W_k(A_1) = \conv \Sigma_k(A_1)$.
We have $\mu\in \Sigma_k(A_1)$.
Then, $\mu = \mu_1 + \cdots + \mu_k$ for eigenvalues $\mu_j$ of $A_1$ corresponding to linearly
independent (unit) eigenvectors $x_1, \ldots, x_k$ of $A_1$. Since $\re \mu_j = \la H_1x_j, x_j\ra \le a_1$ for all $j$,
$\re \mu = n_1a_1$ implies all $x_j\in \cH_1$.
We have ${\rm span}\{x_1, \ldots, x_k\} = \cH_1$. As a linear span of eigenvectors,
 $\cH_1$ is invariant under $A_1$. It is also invariant under
 $G_1 = -i(A_1 - H_1)$. It is indeed a reducing subspace because $G_1$ is self-adjoint.

If $H_1$ has another positive eigenvalue $a_2$, then $\cH_2$ is also a reducing subspace of $H_1$.
This time, we take $k = n_1 + n_2$. The same argument as above yields that $\cH_1\oplus \cH_2$ is a invariant subspace of $A_1$. Consequently, $\cH_2$ is a reducing subspace of $G_1$.
Successively, we get that all $\cH_j$'s are reducing subspaces of $G_1$.

The same can be done for the negative eigenvalues and the corresponding eigenspaces. The only difference is to consider
the left hand supporting line instead. For example, if $b_1$ is the smallest negative eigenvalue of
$B$, then for $k = m_1 = \dim \cK_1$,
$\min\{\re \lambda: \lambda\in W_k(A)\} = m_1b_1$. By considering an extreme point of $\cl W_k(A_1)$ lying on the vertical line $\re \lambda = n_1b_1$, we get that
$\cK_1$, and then similarly, all other $\cK_\ell$'s, are reducing subspaces of $G_1$.
As observed above, $A_1$ is normal.

Similarly, all other $A_2, \ldots, A_m$ are normal. If $A_j = H_j + iG_j$ for self-adjoint
$H_j$ and $G_j$, and if $a + ib$ ($a, b\in\mathbb R$) is an eigenvalue of $A_j$ with corresponding eigenvector $x$, then
$G_j x = ax$ and $H_j x = bx$. Thus we have
$W_k(H_1, G_1, \ldots, H_m, G_m) = \conv \Sigma(H_1, G_1, \ldots, H_m, G_m)$ for all $k$.
For any two operators $B, C\in\{H_1, G_1, \ldots, H_m, G_m\}$, we may apply the argument above to $B + iC$
to conclude that $B$ and $C$ commute. It follows that the operators $A_1, \ldots, A_m$ commute.
\qed

In contrast to the finite dimensional case,
the following example shows that
there are commuting compact self-adjoint operators $H, G$ such that
$W_k(H, G) \equiv W_k(H+iG)$ is not polyhedral for all $k$.

\begin{example} \label{3.2}
Let $A_0 = -I + \diag(e^{i\pi/2}, e^{-i\pi/2},
e^{i\pi/3}, e^{-i\pi/3}, e^{i\pi/4}, e^{-i\pi/4},\dots)$, and
for $k \ge 1$
$A_{k} = \diag(i,-i, -1+i, -1-i)\oplus A_{k-1}$.
Then
$\cl W(A_0)$ is not polyhedral as it has a smooth extreme point at 0.
For $k\ge 1$, $W_{k}(A_k)$ is not polyhedral because the extreme points
$(k-1)i$ and $(1-k)i$ are smooth extreme points.
\end{example}

To determine a condition for
$W_k(\bA)$ to be polyhedral for all $k$,
recall the following notations.
For $\bA = (A_1, \dots, A_m) \in \BH^m$, let
$\bA \oplus {\bf 0}_\ell
= (A_1 \oplus 0_\ell, \dots, A_m\oplus 0_\ell)$,
where $0_\ell$ is the zero matrix in $M_\ell$ if $\ell \in \IN$.
We also write $0_\infty$ as the zero operator on some infinite dimensional Hilbert space.

\begin{theorem} \label{7.3}
Let $\bA = (A_1, \dots, A_m) \in \BH^m$ be
an $m$-tuple of compact operators.
Then for any  $k\in \IN$,
$$\cl(W_k(\bA)) = W_k(\bA \oplus \b0_k) = W_k(\bA \oplus \b0_\infty),$$
and hence
\begin{equation}\label{wk-1}\cl(\conv W_k(\bA))
= \conv W_k(\bA \oplus \b0_k) = \conv W_k(\bA \oplus \b0_\infty).
\end{equation}
Consequently, the following conditions are equivalent.
\begin{itemize}
\item[{\rm (a)}] For every $k \in \IN$,
$\cl W_k(\bA)$ is polyhedral.
\item[{\rm (b)}] For every $k \in \IN$,
$\cl (\conv W_k(\bA))$ is polyhedral.
\item[{\rm (c)}]  The set
$\{A_1, \dots, A_m\}$ is a commuting family of normal operators,
and
$\conv \Sigma_k(\bA \oplus \b0_\infty)$
is polyhedral for every $k \in \IN$,.
\end{itemize}
\end{theorem}

\it Proof. \rm The first assertion follows from Theorem \ref{4.3} and
Theorem \ref{closure}.
One can readily adapt the proof of \cite[Theorem 4.4]{LPW}
to the infinite dimensional operators and show that
for every rank $k$ orthogonal projection
$P = \begin{pmatrix} P_{11} & P_{12} \cr P_{21} & P_{22}\end{pmatrix}$
with $P_{11}, P_{22} \in \BH$,
$P_{11} \oplus P_{22}$ is a convex combination of
orthogonal projections of the form
$Q_1 \oplus Q_2 \in \BH \oplus \BH$, where
$Q_1$ and $Q_2$ are orthogonal projections of rank $\ell$ and $k-\ell$
with $0\le \ell \le k$. Consequently,
$$ W_k(\bA\oplus \b0_\infty)
\subseteq
\conv\bigcup_{0 \le \ell \le k}
(W_\ell(\bA) +  W_{k-\ell}(\b0_\infty))
=\cl (\conv W_k(\bA)).$$
Thus, $\conv W_k(\bA \oplus \b0_k)\subseteq \cl(\conv W_k(\bA))$
and (\ref{wk-1}) follows.

To prove the equivalence of conditions (a) -- (c), we can focus on
self-adjoint operators $A_1, \dots, A_m$  because
we can write $A_\ell = H_\ell + i G_\ell$ for
self-adjoint $H_\ell$ and $G_\ell$ so that

(1) $W_k(\bA)$ is polyhedral if and only if
$W_k(H_1, G_1, \ldots, H_m, G_m)$ is, and

(2) $\{A_1, \dots, A_m\}$ is a commuting family of
normal operators if and only if
$\{H_1, G_1, \ldots, H_m, G_m\}$ is a commuting family
by Fuglede's theorem.

The implication (a) $\Rightarrow$ (b) is clear.
Suppose (b) holds.
For notational simplicity, we let $\bA = \bA \oplus \b0_\infty\in \SH^m$, and
show that the components of $\bA$ commute using an argument similar to those
in the proof of Theorem \ref{3.1}.
Without loss of generality, we shall just show that $A_1A_2 = A_2A_1$.
Now, $\conv W_k(A_1, A_2)$, as the projection
of $\conv W_k(\bA)$ onto $\mathbb R^2$,
is polyhedral.
Since $A_1$ is self-adjoint and compact, it can be written as
$$A_1 = \sum_{j \ge 1} a_j I_{\cH_j} + \sum_{\ell \ge 1}
b_\ell I_{\cK_\ell},$$
where $a_j$'s are the positive, and $b_\ell$'s the negative, eigenvalues
of $A_1$, as in the proof of Theorem \ref{3.1}
We shall consider $A = A_1 + iA_2$ and use \cite[Theorem 1]{MF}
to show that $\cH_j$'s and $\cK_\ell$'s are all reducing subspaces
of $A_2$. It follows that
the null space of $A_1$ is also a reducing subspace of $A_2$, and
consequently, $A_1$ and $A_2$ commute.

If $A_1$ has at least one positive eigenvalue, we start with the largest
eigenvalue $a_1$ and let $k = \dim \cH_1 = n_1$.
For any unit vector $y\in \cH_1^\perp$, let $B = B_1 + iB_2$ be the
compression of
$A = A_1 + iA_2$ onto the subspace $\cH_1\oplus \mbox{span}\{y\}$.
Clearly,
$$\max\{\re \lambda: \lambda\in W_k(B)\} = n_1a_1
= \max\{\re \lambda: \lambda\in W_k(A)\}.$$
Let $\mu$ be an extreme point of $\cl W_k(A)$ lying on the vertical
line $\re \lambda = n_1a_1$.
Then as in the proof of Theorem \ref{3.1}, $\mu\in W_k(A)$.
Suppose $\mu = \sum_{j = 1}^k \la Ax_j, x_j\ra$ for orthonormal vectors
$x_1, \ldots, x_k$. It is easy to see that
all $x_1, \ldots, x_k\in \cH_1$ and so $\mu\in W_k(B)$. As
$W_k(B)\subseteq \cl W_k(A)$, it follows that $\mu$ is a conical
point of $W_k(B)$. By \cite[Theorem 1]{MF},
$\la Bx, y\ra = \la By, x\ra = 0$, or, $\la Ax, y\ra = \la Ay, x\ra = 0$,
for all $x\in \cH_1$.
As this is true for all $y\in \cH_1^\perp$, $\cH_1$ is a reducing subspace
of $A$. It must also be a reducing subspace of $A_2$.

If $A_1$ has another positive eigenvalue $a_2$, then we take $k = n_1 + n_2$ and for any unit vector
$y\in (\cH_1\oplus \cH_2)^\perp$, consider
the compression of $A$ onto $\cH_1\oplus \cH_2\oplus \mbox{span}\{y\}$.
The same argument as above yields that $\cH_1\oplus \cH_2$ is a reducing
subspace of $A_2$. Consequently, $\cH_2$ is a reducing subspace of $A_2$.
Successively, we get that all $\cH_j$'s are reducing subspaces of $A_2$.
Similarly, every $\cK_\ell$ is also reducing for $A$.
Thus, we have $A_1A_2 = A_2A_1$.
Therefore, $\{A_1\oplus 0_\infty, \dots, A_m\oplus 0_\infty\}$, and consequently
$\{A_1 , \dots, A_m\}$, are commuting families.
Condition (c) holds.

Suppose (c) holds.  Since $\{A_1 , \dots, A_m\}$ is a
commuting family of normal operators, by Theorem \ref{3.1},
$W_k(\bA) = \conv \Sigma_k(\bA )$ is convex for every $k \in \IN$.
Now $\{A_1 \oplus 0_\infty, \dots, A_m \oplus 0_\infty\}$ is also a
commuting family of normal operators. Therefore,
$W_k(\bA \oplus \b0_\infty) =\conv \Sigma_k(\bA \oplus \b0_\infty)$
is polyhedral for each $k\in \IN$. By (\ref{wk-1}),
$\cl W_k(\bA) = \cl(\conv W_k(\bA)) = W_k(\bA  \oplus \b0_\infty)$
is also polyhedral.
Hence, condition (a) holds. \qed

For finite rank operators,
the result resembles that of the finite dimensional case.

\begin{proposition} \label{3.4}
Let $\bA = (A_1, \dots, A_m) \in \BH^m$ be
an $m$-tuple of finite rank operators.
  The following conditions are equivalent.
\begin{itemize}
\item[{\rm (a)}]  $\{A_1, \dots, A_m\}$ is a commuting family of
normal operators, equivalently, there is a unitary operator $U$ such that
$U^*A_jU$ is a diagonal operator for each $j$.
\item[{\rm (b)}] $W_k(\bA) = \cl(W_k(\bA))$ is polyhedral for all $k \in \IN$.
\item[{\rm (c)}] There is $k \ge r$ such that
$\conv W_k(\bA)$
 is polyhedral, where  $r$ is the dimension of the sum
of the range spaces of $A_1, \dots, A_m$.
\end{itemize}
\end{proposition}

\it Proof. \rm Suppose (a) holds.
We may assume that $A_1, \dots, A_m$ are real
diagonal operators, each has finite rank.
Then $\Sigma_k(\bA)$ is a finite set. By Theorem \ref{3.1},  $  W_k(\bA) = \conv\Sigma_k(\bA)$ is polyhedral.
Hence $\cl W_k(\bA) = W_k(\bA)$.

The implication (b) $\Rightarrow$ (c)
 is clear.

 Suppose (c) holds.  Then we may assume that
 $A_j = B_j \oplus 0_{\infty}$ with $B_j \in M_r$ for $j = 1, \dots, m$. We may
 further assume that all $B_j$'s are Hermitian. Let $C_j = B_j \oplus 0_{2k-r}$ for
 $j = 1, \dots, m$, $\bB=(B_1,\dots,B_m)$ and $\bC=(C_1,\dots,C_m)$.
By (\ref{wk-1}), we have
$$\cl(\conv W_k(\bA))=\conv W_k(\bB+\b0_k)
\subseteq\conv W_k(\bC)\subseteq\conv W_k(\bA).$$
Therefore, $\conv W_k(\bC)=  \conv W_k(\bA)$ is polyhedral.
For every $1\le i<j\le m$, $W_k(C_i,C_j)=  \conv W_k(C_i,C_j)$ is
a polygon in $\IR^2$. Since $C_i,C_j\in M_{2k}$,
by Theorem \cite[Theorem 3.3]{LPW}, $C_iC_j=C_jC_i$, which implies that $A_iA_j=A_jA_i$.
Hence, (a) holds.
\qed

\section{General operators}
\setcounter{equation}{0}

For any compact convex subset $\cS$ of $\IC^m$, one can construct
diagonal operators $A_1, \dots, A_m$ with
orthonormal basis $\{ e_\mu: \mu  = (\mu_1, \dots, \mu_m) \in \cS\}$
and $A_j e_\mu = (\mu_j/k) e_\mu$ such that
$W_k(I_k\otimes A_1, \dots, I_k \otimes A_m) = \cS$.
(Here and in the following, when we say $A_1, \dots, A_m$ is an $m$-tuple of diagonal
operators, we always assume that they are diagonal operators relative to a fixed
orthonormal basis.)
Thus, there are a lot of flexibility in constructing commuting normal operators
$\bA \in \BH^m$ with prescribed $W_k(\bA)$.
Nevertheless, one can use the geometrical properties
of $W_k(A_1, \dots, A_m)$ to get useful information of $A_1, \dots, A_m$ as
in Proposition \ref{2.1} and Theorem \ref{2.3}, and additional results were obtained for
compact operators in the last section.
However, the situation for general operators are more complicated
as one would expect. For example,
Theorem 6.3 is not true if the operators $A_1, \dots, A_m$ are not compact.
As a counterexample, we exhibit a nonnormal operator whose $k$-numerical range
is a triangle for every $k \in \IN$.

\begin{example} \label{4.1} Let
$A = \begin{pmatrix}0 & 1 \cr 0 & 0 \cr\end{pmatrix}
\oplus [I_\infty \otimes \diag(1, \omega, \omega^2)]$
with $\omega = e^{i2\pi/3}$.
It is not difficult to see that $W_k(A)
 = \conv\{k, k\omega, k\omega^2\}$
is the triangular region with vertices
$k, k\omega, k\omega^2$. However $A$ is not normal because of the first summand.
\end{example}

Despite the above example, we will obtain some results on
$A_1, \dots, A_m$ if $W_k(A_1, \dots, A_m)$ is polyhedral for all $k \in \IN$.

\begin{theorem} \label{4.2}
Let $\bA = (A_1, \dots, A_m)\in \BH^m$. The following conditions are equivalent.
\begin{itemize}
\item[{\rm (a)}]
$W_k(\bA)$ is polyhedral for all $k \in \IN$.
\item[{\rm (b)}] $\conv W_k(\bA)$ is polyhedral for all $k \in \IN$.
\item[{\rm (c)}]
There is a unitary $U\in\BH$ such that $U^*A_jU = D_j \oplus C_j$,
where $D_j = \diag(d_1^{(j)}, d_2^{(j)}, \dots)$ and
$W_k(\bA) = W_k(D_1,\dots, D_m)$ is polyhedral for all
$k\in \IN$.
\end{itemize}
\end{theorem}

\it Proof.  \rm
We may assume that $(A_1, \dots, A_m)
\in \SH^m$.
The implications (c) $\Rightarrow$ (a) $\Rightarrow$ (b) are clear.

Suppose (b) holds.
By Theorem \ref{2.3}, we get for each $k$ an orthonormal set of common
reducing eigenvectors $S_k = \{x_1^{(k)}, \ldots, x_{n_k}^{(k)}\}$ of
$A_1, \ldots, A_m$ such that $W_k(\bA) = W_k(D_1^{(k)}, \dots, D_m^{(k)})$,
where $D_j^{(k)}$ are the compressions of $A_j$ onto ${\rm span}\, S_k$.
We may assume that $S_k\subseteq S_{k + 1}$
by arguing as follows. If $x_1^{(k + 1)}$ is orthogonal to $S_k$, then
$S_k\cup \{x_1^{(k + 1)}\}$ is an orthonormal basis of
${\rm span}\, S_k\cup\{x_1^{(k+1)}\}$.
If $x_1^{(k + 1)}$ is not orthogonal to some of the vectors in $S_k$,
say, $x_1^{(k)}, \ldots, x_\ell^{(k)}$ for $\ell\le n_k$,
then they are all eigenvectors corresponding to the same eigenvalue of
$A_j$ for each $j$.
It is because eigenvectors corresponding to distinct eigenvalues are
orthogonal.
Delete $x_1^{(k + 1)}$ if $x_1^{(k + 1)}\in {\rm span}\,S_k$; otherwise
replace it by a  unit vector $x$ in
${\rm span}\,\{x_1^{(k)}, \ldots, x_\ell^{(k)}, x_1^{(k + 1)}\}$ so that
$S_k\cup \{x\}$ is an orthonormal basis of
${\rm span}\, S_k\cup\{x_1^{(k+1)}\}$.
This $x$ is also a reducing eigenvector of each $A_j$ corresponding to the
same eigenvalue as $x_1^{(k + 1)}$.
Then either $S_k$ or $S_k\cup \{x\}$ is an orthonormal basis of
${\rm span}\, S_k\cup\{x_1^{(k+1)}\}$.
Going through all other vectors in $S_{k + 1}$,
we get an orthonormal set $S_{k + 1}'$ containing $S_k$ such that
${\rm span}\, S_{k + 1}' = {\rm span}\, S_k\cup S_{k + 1}$.
Note that $S_{k + 1}'$ contains all vectors in $S_{k + 1}$ that are
orthogonal to $S_k$, and for a vector that is not
orthogonal to $S_k$, it is either deleted or replaced.
A vector that is deleted can be regarded as replaced by an eigenvector in
$S_k$ corresponding to the same eigenvalue.
With the previous notations, if $x_1^{(k + 1)}, \ldots, x_r^{(k + 1)}$
belong to ${\rm span}\,\{x_1^{(k)}, \ldots, x_\ell^{(k)}\}$ so that they
are deleted, then $r\le \ell$. There are always enough vectors for
replacement.
Renaming $S_{k + 1}'$ as $S_{k + 1}$, then we still have
$W_{k + 1}(\bA) = W_{k + 1}(D_1^{(k+1)}, \dots, D_m^{(k+1)})$.
Thus we may assume that the sets $S_k$ are increasing. It follows that
$S = \bigcup_k S_k$ is an orthonormal set
of common reducing eigenvectors of $A_1, \ldots, A_m$.
If $D_j$ is the compression of $A_j$ onto ${\rm span}\, S$,
then $W_k(\bA) = W_k(D_1, \dots, D_m)$ is polyhedral  for all $k\in \IN$.
Hence, condition (c) holds.
\qed

\begin{theorem} \label{7.3b}
Suppose $\bA \in \BH^m$. The following are equivalent.
\begin{itemize}
\item[{\rm (a)}]  $\cl W_k(\bA)$ is polyhedral
for every $k\in \IN$.
\item[{\rm (b)}]  $\cl(\conv W_k(\bA))$ is polyhedral
for every $k\in \IN$.
\item[{\rm (c)}]
There exists an $m$-tuple of diagonal operators
$\bD$ of the form $(R_1 \oplus S_1, \dots, R_m\oplus S_m)$,
where $\bR = (R_1, \dots, R_m)$
is the compression of $\bA$ on a separable reducing subspace of $\cH$,
and $\bS = (S_1, \dots, S_m)$ is an $m$-tuple of diagonal operators
whose diagonal entries are extreme points of
$W_{\rm ess}(\bA)$, such that
$$\cl W_k(\bA) = W_k(\bD) \quad \hbox{ is polyhedral
for every } \ k \in \IN.$$
\end{itemize}
\end{theorem}

\it Proof. \rm
We may assume that $(A_1, \dots, A_m)\in \SH^m$.
The implications (c) $\Rightarrow$ (a) $\Rightarrow$ (b) are clear.

Suppose (b) holds. For each $k\in \IN$, let $\bv_1^k,\dots,\bv_{n_k}^k$
be the vertices of $\cl (\conv W_k(\bA))$.  For each vertex $\bv_j^k$, $1\le j\le n_k$,
by Theorem \ref{closure}, there exist $0\le \ell\le k$, $\bw^k_j\in W_{\ell}(\bA)$ and
$\bd^k_j\in W_{\rm ess}(\bA)$ such that
 $\bv_j^k=\bw^k_j+(k-\ell)\bd^k_j$.  Let
$$\bS=\oplus_{k\ge 1}\oplus_{1\le j\le n_k}\bd^k_jI_{\infty}.$$
Then for every  $k\in \IN$,
$\conv W_k(\bA\oplus \bS)=\cl W_k(\bA) $
is polyhedral for each $k \in \IN$. By Theorem \ref{4.2},
$\bA \oplus \bS$ is unitarily similar to
$\bD \oplus \bC$ for an $m$-tuple of diagonal operators $\bD$ such that
$W_k(\bA \oplus \bS) = W_k(\bD)$ for each
$k \in \IN$.
Let $\cS = \{(x_j, y_j)\}$ be an orthonormal set of common reducing eigenvectors
of $\bA \oplus  \bS$ corresponding to the diagonal entries of $\bD$.
Then $\cS$ is at most countably infinite.
If $y_j = 0$, then $x_j$   is a common reducing eigenvector of $\bA$.
Let $\cS_1 = \{x_j:(x_j, 0)\in \cS\}$. Then $\cS_1$ is orthonormal.
Extending $\cS_1$ to an orthonormal basis of $\cH$, we have a unitary $U$
in $\BH$ and $m$-tuple of $\bR = (R_1, \dots, R_m)$ of
diagonal operators acting on a separable subspace of $\cH$
such that
$U^*A_jU = R_j \oplus B_j$ for $j = 1, \dots, m$.  Note that  the common
eigenvalue $\bd_j$ corresponding to each $x_j\in \cS_1$ is included in $\bR$.
Also, it is easy to show that it suffices to use the extreme points
of $W_{\ess}(\bA)$ to construct the diagonal operator $\bD$.
Then $\bR$ and $\bS$ satisfy condition (c). \qed

One may wonder whether $\bR$ or $\bS$ in
Theorem 7.3 can be replaced by finite dimensional operators.
By Proposition \ref{W(A)/k},
for any $\bA \in \BH^m$
$$\frac{1}{k+1}W_{k+1}(\bA) \subseteq
\cl (\conv \frac{1}{k}W_{k}(\bA)) \quad \hbox{ if } k \in \IN.$$
Hence,
$$\cl (\conv \frac{1}{k}W_{k}(\bA)) \rightarrow W_{\rm ess}(\bA)\quad
\hbox{as}\quad k \rightarrow \infty.$$
One may ask whether $\cl(W_k(\bA))$ is polyhedral for all $k \in \IN$
will imply that $W_{\ess}(\bA)$ is polyhedral, and hence the distinct
diagonal entries of $\bS$ form a finite set.
The following example shows that the
answers for the above questions are both negative.

\begin{example}
Let $D_1 = \diag(1, e^{i\pi/2}, e^{-i\pi/2}, e^{i\pi/3}, e^{-i\pi/3}, \dots)$,
$D_2 = (1)\oplus \diag(1 + 1, 1 + 1/2, 1 + 1/3, \dots)$, $A_1 = D_1\otimes I_\infty$ and
$A = A_1\oplus D_2$.
Then $W_{\rm ess}(A) =
\conv\{1, e^{i\pi/2}, e^{-i\pi/2}, e^{i\pi/3}, e^{-i\pi/3}, \dots\}$
is not polyhedral. On the other hand, it is easy to see that
$W_1(A) = \conv (W_1(A_1)\cup W_1(D_2))$ is polyhedral. For $k \in \IN$, we have
$\frac 1k W_k(A_1) = W_{\ess}(A) = \conv\{1, e^{i\pi/2}, e^{-i\pi/2}, e^{i\pi/3}, e^{-i\pi/3}, \dots\}$,
and $\frac 1k W_k(D_2)$ is a line segment on the real axis of the form $[1, 1 + r]$
for some $r > 0$. One may use
\cite[Theorem 4.4]{LPW}
and mathematical induction
to show that $\frac 1k W_k(A)$ has only finitely many extreme points.
\end{example}

Theorems \ref{4.2} and \ref{4.3} characterize $\bA$ such that  $\cS_k(\bA)$
is polyhedral for all $k \in \IN$, where $\cS_k(\bA)$ is $W_k(\bA)$,
$\conv W_k(\bA)$, $\cl W_k(\bA)$, or  $\cl(\conv W_k(\bA)$.
In all cases, the condition depends on the
behavior of the compression of $A_1, \dots, A_m$ on a certain
separable subspace of $\cH$.
Even if $\cH$ is separable, Example \ref{4.1} shows that
it is not easy to determine whether $\{A_1, \dots, A_m\}$
is a commuting normal family using the geometrical structure of
$W_k(\bA)$ for all $k \in \IN$. We are only able to obtain a
result for some special families of $\{A_1, \dots, A_m\}$,
which can be deduced from the following.

\begin{proposition} \label{4.4}
Let $A \in \BH$ have a finite spectrum.
Then $A$ is normal if any one of the following conditions holds.
\begin{itemize}
\item[{\rm (a)}]
For each $\lambda \in \sigma(A)$ there is  $k\in \IN$ such that
$W_k(A)$ is polyhedral with a conical point
$\lambda + \sum_{j=2}^k \mu_j \in \Sigma_k(A)$.
\item[{\rm (b)}]
For each $\lambda \in \sigma(A)$ there is $k \in \IN$ such that
$\cl W_k(A)$ is polyhedral with a conical point
$\lambda + \sum_{j=2}^\ell \mu_j + (k-\ell)\xi$, where $1\le\ell \le k$,
$\lambda+\sum_{j=2}^\ell \mu_j \in \Sigma_\ell(A)$ and $\xi \in W_{\ess}(A)$.
\end{itemize}
\end{proposition}

\it Proof. \rm Suppose $A \in \BH$ has a finite spectrum
and (a) holds. Let $\lambda \in \sigma(A)$ and
$W_k(A)$ be polyhedral with a vertex of the form
$\sum_{j=1}^k \mu_j \in \Sigma_k(A)$, where
$\mu_1 = \lambda$. We will show that $A$ is unitarily
similar to $\lambda I \oplus B$ such that $\lambda \notin \sigma(B)$.
By Theorem \ref{2.3}, $A$ is unitarily similar to $D \oplus B$ for a
diagonal matrix $D\in M_r$ such that $W(A) = W(D)$. In particular, we may
assume that $\mu_1, \dots, \mu_m$ are diagonal entries of $D$.
If $\lambda \notin \sigma(B)$, then we are done.
Otherwise, assume that $\lambda \in \sigma(B)$. We may assume that
the first column (of the operator matrix) of $B$ has the form
$(\lambda, 0, \dots)^t$. We claim that the first
row $B$ is of the form $(\lambda, 0, \dots)$.
If not, we may assume that the $(1,2)$ entry of $B$ is $\gamma \ne 0$.
Consider the leading principal submatrix $T\in M_{r+2}$
of $D \oplus B$. Then the $(r+1,r+2)$ entry of $T$ is $\gamma$.
Note that
$$W_k(D) \subseteq W_k(T) \subseteq W_k(A) = W_k(D).$$
We can find $k$ diagonal entries of $T$ including the $(r+1,r+1)$
entry, which is $\lambda$, summing up to $\sum_{j=1}^k \mu_j$
is a vertex of $W_k(T)$.  By Theorem \ref{5.4} (a),
the $(r+1,r+2)$ entry of $T$ is zero, which is a contradiction.
So, if $\lambda \in \sigma(B)$, then $\sigma$ is a reducing eigenvalue if $B$.
As a result, $B= \lambda I \oplus B_1$ such that
$\lambda \notin \sigma(B_1)$. Repeating this argument for every
element in $\sigma(A)$, we see that $\cH$ is a direct sum of the eigenspaces
of $A$. Thus, $A$ is normal.

Next, suppose (b) holds.
Let $\lambda \in \sigma(A)$, and $\cl W_k(A)$ is polyhedral with a conical
point $\lambda + \sum_{j=2}^\ell \mu_j + (k-\ell) \xi$.
By Theorem \ref{2.4} and its proof, we see that
$\cl W_\ell(A)$ is polyhedral with $\lambda + \sum_{j=2}^\ell \mu_j$
as a conical point. We may use the argument in the preceding
paragraph to show that $\lambda$ is a reducing eigenvalue of $A$.
Since this is true for all $\lambda \in \sigma(A)$, we see that $A$ is normal.
\qed

In Proposition \ref{4.4} (a) and (b), one cannot just require $W_k(A)$ to be polyhedral
for $k \le \ell$ for a finite value $\ell$ using
the information of $\sigma(A)$.
For example, suppose $\sigma(A) = \{1, w, w^2\}$ with $w = e^{i2\pi/3}$.
For any  $\ell\in \IN$, we can let
$A = [I_\ell \otimes \diag(1, w, w^2)] \oplus
\begin{pmatrix} 0 & 1 \cr 0 & 0\cr\end{pmatrix}$
to be the non-normal operator
such that $W_k(A) = \conv\{k, kw, kw^2\}$
is a polygon for all $k = 1, \dots, \ell$, but
$W_{k+1}(A)$ is not.

\begin{corollary} \label{4.5}
Let $\bA = (A_1,\dots, A_m) \in \BH^m$ be an
$m$-tuple of  self-adjoint operators. Then
$\{A_1, \dots, A_m\}$ is a commuting family if
for all choices of  $1 \le p < q \le m$,
 $A = A_p + iA_q$ has a finite spectrum  satisfying  conditions {\rm (a)} or {\rm  (b)}
 of Proposition {\rm \ref{4.4}}.
\end{corollary}

The following example shows that Corollary
\ref{4.5} will fail if one only assumes that
$A_1, A_2$ are self-adjoint operators, where
$A_1, A_2$ have finite spectrum, and for any eigenvalues
$\lambda_j \in \sigma(A_j)$ for $j = 1,2$,
there is $k \in \IN$ such that $W_k(A_1+iA_2)$ is
polyhedral with  a vertex of the form
$(\lambda_1 + i \lambda_2) + \sum_{j=2}^k \mu_j \in \Sigma_k(A_1+iA_2)$.

\begin{example} \label{7.7}
Let $A_1= \begin{pmatrix}1 & 0 \cr 0& -1 \cr\end{pmatrix}
\oplus (I_\infty \otimes D_1)$ and
$A_2= \begin{pmatrix}0 & 1\cr 1& 0 \cr\end{pmatrix}
\oplus (I_\infty \otimes D_2)$ with $D_1 = \diag(1, 1,-1,-1)$,
$D_2 = \diag(1,-1,1,-1)$.
Then $A_1A_2 \ne A_2A_1$,
$\sigma(A_1)=\sigma(A_2)=\{1,-1\}$. Let $A=A_1+iA_2$. Then for every $k\ge 1$,
$W_k(A)$ is the square with vertices $\{k(1+i),k(1-i),k(-1+i),k(-1-i)\}$. Note that $0\in \sigma(A_1+iA_2)$ but a point   $\sum_{j=1}^k\mu_j$ is a vertex of $W_k(A)$ if and only if
$\mu_1=\cdots=\mu_k\in  \{ 1+i , 1-i , -1+i ,k -1-i \}$.
\end{example}

\section{Open problems and related results}
\setcounter{equation}{0}
There are many interesting problems and related results in connection to our study.
We list some of them in the following.
\begin{itemize}
\item In connection to Proposition \ref{W(A)/k}, it is interesting to prove or disprove  that
$W_{k+1}(\bA)/(k+1) \subseteq W_k(\bA)/k$.
\item In connection to Theorem \ref{4.3} and Theorem \ref{k-closed}, it is interesting to
prove or disprove that
$W_{k}(\bA)$ is closed whenever $W_{k+1}(\bA)$ is closed.

\item In connection to Example \ref{4.6}, it is interesting to construct an
example $\bA \in \BH^n$ such that $\conv W_k(\bA)$ is closed, but $W_k(\bA)$ is not
for $k > 1$.
\end{itemize}

For any real vector $c = (c_1, \dots, c_k)$ one can define the joint
$c$-numerical range of $\bA\in \BH^m$
as the collection of $(a_1, \dots, a_m) \in \IC^m$
such that $a_j = \sum_{j=1}^k c_j \la A_j x_j, x_j\ra$ for an orthonormal set
$\{x_1,\dots, x_k\} \subseteq \cH$.
One may see \cite{LP0,LST} for some background
of the $c$-numerical range for a single operator.
Some of our results can be readily extended to
the setting of the joint $c$-numerical ranges. For example,
if $c = (c_1, \dots, c_k)$ has distinct entries, then $W_c(\bA)$ has a conical
boundary point $(a_1, \dots, a_m)$ will imply that there is a unitary $U \in\BH$ such that
$U^*A_jU = D_j \oplus B_j$ for $j = 1, \dots, m$,
where $D_j \in M_k$ is a diagonal matrix satisfying
$(\tr(CD_1), \dots, \tr(CD_m))=(a_1, \dots, a_m) $ with $C = \diag(c_1, \dots, c_k)$.
It would be interesting to extend other results to $W_c(\bA)$.
For example, it was shown in \cite{LST}
for $c = (c_1, \dots, c_n)$ with entries arranged and descending
order and $c_k > c_{k+1}$, where $|n/2-k| \le 1$, then
 $A \in \bM_n$ satisfies $W_c(A)$ is polyhedral if and only if
 $A$ is normal. One can  readily show that
 $W_c(\bA)$ is polyhedral for $(c_1, \dots, c_n)$
 satisfying the above condition and $\bA = (A_1, \dots, A_m) \in M_n^m$
 if and only if $\{A_1, \dots, A_m\}$
 is a commuting family of normal matrices; it follows that
 $W_d(\bA)$ is polyhedral for all $d = (d_1, \dots, d_n) \in \IR^n$.
It would be interesting to extend the results
to $W_c(\bA)$ for $\bA\in \BH^m$ or $\bA \in \KH^m$.

In \cite{Cho2}, the authors studied  special boundary points of $W(\bA)$ such as
conical points, bare points, or points with some extreme norm properties.
We have extended their result for conical points to $W_k(\bA)$. It is interesting to
extend other results in the paper to $W_k(\bA)$ or $W_c(\bA)$.

\bigskip\noindent
{\Large\bf No conflict of interest statement}

 On behalf of all authors, the corresponding author states that there is no conflict of interest.

\section*{Acknowledgement}

Li is an affiliate member of the Institute for Quantum
Computing, University of Waterloo; his
research was partially supported by the
Simons Foundation Grant 851334. The authors would like to thank the
anonymous referee for some helpful comments.

\bigskip\noindent
[Chan]
 {jtchan@hku.hk}

\medskip\noindent
[Li] {Department of Mathematics, The College of William
\& Mary, Williamsburg, VA 13185, USA.}  {ckli@math.wm.edu}

\medskip\noindent
[Poon] {Department of Mathematics, Iowa State University,
Ames, IA 50011, USA.}
{ytpoon@iastate.edu}
\end{document}